\documentclass[final]{siamart0516}

\usepackage{siampaper}
\usepackage{mathabx}
\usepackage[final]{changes}

\usepackage[T1]{fontenc}

\renewcommand{\year}{2020}     % If you don't want the current year.
\newcommand{\cahiernumber}{72}  % Insert your Cahier du GERAD number.

\newsiamthm{principle}{Principle}
\Crefname{algorithm}{Algorithm}{Algorithms}
\Crefname{ALC@unique}{Line}{Lines}
\Crefname{assumption}{Assumption}{Assumptions}
\newcommand{\diag}{\mathop{\mathrm{diag}}}
\newcommand{\rank}{\mathop{\mathrm{rank}}}

\algnewcommand{\IfOneLine}[1]{\State\algorithmicif\ {#1},}

\title{
  Constraint-Preconditioned Krylov Solvers for Regularized Saddle-Point Systems
}

\author{
  Daniela di Serafino%
  \thanks{%
    Department of Mathematics and Applications "R. Caccioppoli",
    University of Naples Federico II, Naples, Italy.
    E-mail: \mailto{daniela.diserafino@unina.it}.
    Research partially supported by GERAD during a visit of this author in 2017,
    and by Gruppo Nazionale per il Calcolo Scientifico -- Istituto Nazionale di Alta Matematica (GNCS--INdAM), Italy.
  }
  \and
  Dominique Orban%
  \thanks{%
    GERAD and Department of Mathematics and Industrial Engineering,
    \'Ecole Polytechnique, Montr\'eal, QC, Canada.
    E-mail: \mailto{dominique.orban@gerad.ca}.
    Research partially supported by an NSERC Discovery Grant.
  }
}
\date{\today}

\begin{document}

\maketitle

\thispagestyle{firstpage}
\pagestyle{myheadings}

\begin{abstract}
We consider the iterative solution of regularized saddle-point systems.
When the leading block is symmetric and positive semi-definite on an appropriate subspace,
\cite{dollar-gould-schilders-wathen-2006} describe how to apply the conjugate gradient (CG) method coupled
with a constraint preconditioner, a choice that has proved to be effective in optimization applications.
We investigate the design of constraint-preconditioned variants of other Krylov methods for regularized systems
by focusing on the underlying basis-generation process. We build upon principles laid out by~\cite{gould-orban-rees-2014}
to provide general guidelines that allow us to specialize any Krylov method to regularized saddle-point systems.
In particular, we obtain constraint-preconditioned variants of Lanczos and Arnoldi-based methods, including the
Lanczos version of CG, MINRES, SYMMLQ, GMRES($\ell$) and DQGMRES.
We also provide MATLAB implementations in hopes that they are useful as a basis for the development of more
sophisticated software. Finally, we illustrate the numerical behavior of constraint-preconditioned Krylov solvers
using symmetric and nonsymmetric systems arising from constrained optimization.
% and fluid-flow simulation.
\end{abstract}

\begin{keywords}
Regularized saddle-point systems, constraint preconditioners, Lanczos and Arnoldi procedures, Krylov solvers.
\end{keywords}

\begin{AMS}
65F08, 65F10, 65F50, 90C20.
\end{AMS}

\smallskip

% \listoftodos\relax        % Remove when finished

\section{Introduction\label{sec:introduction}}

We consider the iterative solution of the regularized saddle-point system
\begin{equation}
  \label{eq:rsp}
  \begin{bmatrix}
    A & \phantom{-}B^T \\
    B & -C\phantom{^T}
  \end{bmatrix}
  \begin{bmatrix}
    x \\ y
  \end{bmatrix}
  =
  \begin{bmatrix}
    b \\ 0
  \end{bmatrix},
\end{equation}
where $A \in \R^{n \times n}$ may be nonsymmetric, $C \in \R^{m \times m}$ is nonzero and symmetric, and $B \in \R^{m \times n}$.
We denote $K$ the matrix of~\eqref{eq:rsp}.
There is no loss of generality in assuming that the last $m$ entries of the right-hand side of~\eqref{eq:rsp} are zero, as discussed later.

A constraint preconditioner for~\eqref{eq:rsp} has the form
\begin{equation}
   \label{eq:cp}
   P =
   \begin{bmatrix}
   G & \phantom{-}B^T \\
   B & -C\phantom{^T}
   \end{bmatrix},
\end{equation}
where $G$ is an approximation to $A$ such that~\eqref{eq:cp} is nonsingular.
When $A$ is symmetric and has appropriate additional properties, a constraint preconditioner allows the application of CG even
though $K$ and $P$ are indefinite \citep{dollar-gould-schilders-wathen-2006}.
% \citep{gould-hribar-nocedal-2001,dapuzzo-desimone-diserafino-2010}.

We are interested in the design of constraint-preconditioned versions of additional Krylov methods for~\eqref{eq:rsp},
including methods that can be used when $A$ is nonsymmetric. We extend the work of~\cite{gould-orban-rees-2014} on projected
and constraint-preconditioned Krylov methods for saddle-point systems with $C=0$ by exploiting a suitable reformulation of~\eqref{eq:rsp}
suggested by \cite{dollar-gould-schilders-wathen-2006}. We develop constraint-preconditioned variants of the Lanczos and Arnoldi basis-generation
processes, and use them to derive variants of Krylov solvers based on those processes.
More generally, we provide guidelines that can be also exploited to obtain constraint-preconditioned versions of other Krylov methods
not considered in this paper. Finally, we distribute MATLAB implementations of the constraint-preconditioned methods discussed
here as templates for the development of more sophisticated numerical software.

%
%\citep{luksan-vlcek-1998,
%keller-gould-wathen-2000,
%perugia-simoncini-2000,
%durazzi-ruggiero-2003,
%bergamaschi-gondzio-zilli-2004,
%dollar-gould-schilders-wathen-2006,
%cafieri-dapuzzo-desimone-diserafino-2007a,
%bergamaschi-gondzio-venturin-zilli-2007, dollar-2007,
%forsgren-gill-griffin-2007,
%sesana-simoncini-2013,
%gould-orban-rees-2014,
%bellavia-desimone-diserafino-morini-2015,
%bellavia-desimone-diserafino-morini-2016,
%fisher-gratton-gurol-tremolet-vasseur-2016,
%bergamaschi-desimone-diserafino-martinez-2017}.
%
% \citep{orban-2015}?? Not a constraint preconditioner (like work by Scott and Tuma).
%

% In this case, a common choice is $G = \diag(A)$, provided that all the diagonal entries of $A$
% are positive. The matrix $G$ can be also implicitly defined by considering a factorization
% $P = N R \, N^T$, with specially chosen factors $N$ and  $R$ \citep{dollar-gould-schilders-wathen-2006}.
% Furthermore, the preconditioned matrix $P^{-1} M$ has an eigenvalue at~1 with multiplicity (at least)
% $2m-p$, and the remaining eigenvalues are defined by a generalized eigenvalue problem
% with real solutions~\citep{dollar-2007}. Roughly speaking, the better $G$ approximate $A$
% the more clustered around~1 are the remaining eigenvalues (see~\cite{dollar-2007} for details).

Systems of type~\eqref{eq:rsp} arise in interior-point methods for constrained optimization in the presence of inequality constraints
or when regularization is used \citep{benzi-golub-liesen-2005,dapuzzo-desimone-diserafino-2010,friedlander-orban-2012}.
They also appear in Lagrangian approaches for variational problems with equality constraints when the constraints are relaxed
or a penalty term is applied \citep{pestana-wathen-2015}. In the above cases, $A$ is usually symmetric, but may also be
nonsymmetric---see~\Cref{sec:experiments}, and often has additional properties,
e.g., it accounts for local convexity of the optimization problem. Regularized saddle-point systems with nonsymmetric $A$ arise
also from the stabilized finite-element discretization of Oseen problems obtained by linearization, through Picard's method,
of the steady-state Navier-Stokes equations governing the flow of a Newtonian incompressible viscous fluid \citep{benzi-golub-liesen-2005}.

Constraint preconditioners have widely demonstrated their effectiveness
on saddle-point systems, especially when the leading block is symmetric and enjoys
additional properties, such as being positive definite; much work has been
carried out to develop, analyze and approximate constraint preconditioners in this
case, see, e.g., \citep{benzi-golub-liesen-2005,dapuzzo-desimone-diserafino-2010,gould-orban-rees-2014,
desimone-diserafino-morini-2018} and the references therein.
%\citep{benzi-golub-liesen-2005,dapuzzo-desimone-diserafino-2010,gould-orban-rees-2014,
%bellavia-desimone-diserafino-morini-2015, bellavia-desimone-diserafino-morini-2016,
%fisher-gratton-gurol-tremolet-vasseur-2016, bergamaschi-desimone-diserafino-martinez-2017,
%desimone-diserafino-morini-2018}.

The rest of this paper is organized as follows.
\Cref{sec:preliminaries} provides preliminary results used in the sequel.
In \Cref{sec:cp-lanczos}, we describe the constraint-preconditioned Lanczos process and,
in \Cref{sec:cp-l-methods}, we present variants of Krylov solvers based on it.
In \Cref{sec:cp-arnoldi-and-methods}, we describe the constraint-preconditioned Arnoldi process and associated Krylov methods.
In \Cref{sec:implementation}, we discuss implementation issues and provide details on the MATLAB codes.
In \Cref{sec:experiments}, we illustrate the numerical behavior of some constraint-preconditioned solvers on regularized saddle-point systems,
with symmetric and nonsymmetric matrices, from constrained optimization.
% and fluid-flow simulation.
We conclude in \Cref{sec:discussion}.

\subsubsection*{Notation}

Uppercase Latin letters ($A$, $B$, $\ldots$), lowercase Latin letters ($a$, $b$, $\ldots$),
and lowercase Greek letters ($\alpha$, $\beta$, $\ldots$) denote matrices, vectors and scalars, respectively.
The Euclidean norm is denoted $\| \cdot \|$. If \(S = S^T\) is a positive definite matrix, the \(S\)-norm is defined as \(\|u\|_S^2 = u^T S \, u\).
All vectors are column vectors. For any vector $v$, $\diag(v)$ is the diagonal matrix with diagonal entries equal to the entries of $v$.
For brevity, we use the MATLAB-like notation $[v \,;\, w]$ to represent the vector $[v^T \; w^T]^T$.

\section{Preliminaries\label{sec:preliminaries}}

We assume throughout that $K$ is nonsingular, which implies
\begin{equation}
  \label{eq:rsp-nonsingular}
  \Null(A) \cap \Null(B) = \{0\}
  \quad \text{and} \quad
  \Null(B^T) \cap \Null(C) = \{0\}.
\end{equation}
In general the converse is not true.
A counterexample consists in taking
$$
    A = \begin{bmatrix}
       \phantom{-}1 & -1                   & 0 \\
       \phantom{-}0 & \phantom{-}0 & 0 \\
       \phantom{-}1 & \phantom{-}0 & 1
    \end{bmatrix}, \quad
   B = \begin{bmatrix}
      1 & 0 & 0 \\
      0 & 0 & 1
   \end{bmatrix}, \quad
   C = \begin{bmatrix}
      1 & 0 \\
      0 & 1
   \end{bmatrix} .
$$
\cite{benzi-golub-liesen-2005} and \cite{dapuzzo-desimone-diserafino-2010} give additional conditions that guarantee nonsingularity of $K$.
Note however that we do not require \(B\) to have full rank or \(C\) to be positive (semi-)definite.

In order to develop constraint-preconditioned Krylov methods for~\eqref{eq:rsp}, we specialize the basis-generation processes underlying those methods.
We focus on the \cite{lanczos-1950} and \cite{arnoldi-1951} processes, which compute orthonormal bases of Krylov spaces associated with symmetric and general matrices,  respectively.
For reference, the preconditioned Lanczos process is stated as \Cref{alg:prec-lanczos} in \Cref{sec:processes}.
The standard Lanczos process follows by setting the preconditioner to the identity.
It is straighforward to apply our arguments to the \cite{lanczos-1950} biorthogonalization process and its transpose-free variants \citep{brezinski-redivozaglia-1998, chan-depillis-vandervorst-1998}.
We implicitly assume that $A = A^T$ when considering the Lanczos process.

Following \cite{dollar-gould-schilders-wathen-2006}, we reformulate~\eqref{eq:rsp} as follows.
Assume that $\rank ( C ) = p$ and $C$ has been decomposed as\footnote{Note that~\eqref{eq:factc} will be only used for the purpose of deriving computational processes and need not be computed in practice.}
\begin{equation}
\label{eq:factc}
   C = E F E^T\!,
\end{equation}
where $F \in \R^{p \times p}$ is symmetric and nonsingular and $E \in \R^{m \times p}$.
Then, by using the auxiliary variable
\begin{equation}
\label{eq:w}
w = - F E^T y,
\end{equation}
equation~\eqref{eq:rsp} may be written
\begin{equation}
   \label{eq:rsp-2}
   \begin{bmatrix}
      A &                              & B^T \\
         & F^{-1}                    & E^T \\
      B & E\phantom{^{-1}}  &
   \end{bmatrix}
   \begin{bmatrix}
      x \\ w \\ y
   \end{bmatrix}
   =
   \begin{bmatrix}
      b \\ 0 \\ 0
   \end{bmatrix},
\end{equation}
which has a standard symmetric saddle-point form
\begin{equation}
  \label{eq:def-3x3}
  \begin{bmatrix}
       M & N^T
    \\ N &
  \end{bmatrix}
  \begin{bmatrix}
    g \\ y
  \end{bmatrix}
  =
  \begin{bmatrix}
    b_0 \\ 0
  \end{bmatrix},
  \quad
  M =
  \begin{bmatrix}
    A & \\
      & F^{-1}
  \end{bmatrix},
  \
  N =
  \begin{bmatrix}
    B & E
  \end{bmatrix},
  \
  g =
  \begin{bmatrix}
    x \\ w
  \end{bmatrix},
  \
  b_0 =
  \begin{bmatrix}
    b \\ 0
  \end{bmatrix}.
\end{equation}
The principles laid out by \cite{gould-orban-rees-2014} may now be applied to~\eqref{eq:rsp-2}.

Note that~\eqref{eq:rsp-2} is nonsingular if and only if~\eqref{eq:rsp} is nonsingular, and therefore $N$ must have full rank.
Because $g \in \Null(N)$, there exists $\widehat d \in \R^{n+p-m}$ such that
\begin{equation}
\label{eq:nullspace-basis}
    g =
    \begin{bmatrix}
      x \\ w
    \end{bmatrix}
    = Z \, \widehat d =
    \begin{bmatrix}
      Z_1 \\ Z_2
    \end{bmatrix} \widehat d ,
\end{equation}
where the columns of $Z$ form a basis of $\Null(N)$.
%
% In order to derive constraint-preconditioned Lanczos and Arnoldi procedures, we start
% from the formulation of projected versions of them, which in turn are obtained from
% the corresponding versions restricted to $\Null([B \;\, E])$, as briefly shown next.
%
The restriction of~\eqref{eq:rsp-2} to $\Null(N)$ is
\begin{equation}
\label{eq:reduced-sys}
    \widehat{M} \, \widehat{x} = \widehat{b},
\end{equation}
where
\begin{subequations}
  \label{eq:reduced-sys-2}
  \begin{align}
    \widehat{M} & = Z^T M Z = Z_1^T A Z_1 +Z_2^T F^{-1} Z_2,
    \\
    \begin{bmatrix}
    x \\ w
    \end{bmatrix}
    & =
    \begin{bmatrix}
      Z_1 \\ Z_2
    \end{bmatrix}
    \widehat{x},
    \quad
    \widehat{b} =
    \begin{bmatrix}
      Z_1^T & Z_2^T
    \end{bmatrix}
    \begin{bmatrix}
      b \\ 0
    \end{bmatrix}
    = Z_1^T b .
  \end{align}
    %
    % \begin{array}{cc}
    % \displaystyle
    %
    % \begin{bmatrix}
    %   Z_1^T & Z_2^T
    % \end{bmatrix}
    % \begin{bmatrix}
    %   A & \\
    %     & F^{-1}
    % \end{bmatrix}
    % \begin{bmatrix}
    %   Z_1 \\ Z_2
    % \end{bmatrix}
    % \\[4mm]
    % \displaystyle
    %
    % \end{array}
\end{subequations}
% Henceforth, $\widehat{M}$ is assumed to be symmetric when the Lanczos process is analyzed,
% while $\widehat{M}$ may be general when the Arnoldi one is considered.

In a Krylov method for~\eqref{eq:reduced-sys}, it is appropriate to use a preconditioner of the form
\begin{equation}
    \label{eq:reduced-prec}
    \widehat{P} = Z_1^T G Z_1 +Z_2^T F^{-1} Z_2.
    %\begin{bmatrix}
    %  Z_1^T & Z_2^T
    %\end{bmatrix}
    %\begin{bmatrix}
    %  G & \\
    %   & F^{-1}
    %\end{bmatrix}
    %\begin{bmatrix}
    %  Z_1 \\ Z_2
    %\end{bmatrix}.
\end{equation}
If $G$ is suitable, the preconditioned method can be reformulated entirely in terms of full space
quantities \citep{gould-hribar-nocedal-2001,dollar-gould-schilders-wathen-2006,gould-orban-rees-2014}.
Following \citep[Assumption~2.2]{gould-orban-rees-2014}, we require the following assumption.
\begin{assumption}
  \label{assum:Phat-spd}
  % The preconditioner $\widehat{P}$ is symmetric and positive definite.
  The matrix
  $$
  \begin{bmatrix}
    G & \\
      & F^{-1}
  \end{bmatrix}
  $$
  is symmetric and positive definite on $\Null(N)$.
\end{assumption}
% Note that $\widehat P$ is symmetric if $G$ is symmetric.
% Positive definiteness of $\widehat P$ can be verified without computing $Z_1$ and $Z_2$---e.g.,
% \cite[Theorem~2.1]{dollar-gould-schilders-wathen-2006} and \cite[p.~290]{dapuzzo-desimone-diserafino-2010}.
A consequence of \Cref{assum:Phat-spd} is that~\eqref{eq:reduced-prec} is symmetric and positive definite.

We enforce \Cref{assum:Phat-spd} throughout this paper to guarantee that Krylov methods for~\eqref{eq:reduced-sys} give rise to corresponding full-space methods for~\eqref{eq:rsp}.
However, at least in principle, \Cref{assum:Phat-spd} is not always necessary, e.g., in Krylov methods based on the Arnoldi process.
% On the other hand, a positive definite $\widehat P$ can be effective even if $\widehat M$ is nonsymmetric, as shown in \Cref{sec:experiments}.

The application of the preconditioner $\widehat{P}$, i.e., $\widecheck u = \widehat{P}^{-1} \widehat u$, can be written as
% \smarttodo{I wish we could find a better notation than those hats}
\begin{equation}
  \label{eq:apply-PG}
   \begin{bmatrix}
   \bar{u}_x \\ \bar{u}_w
   \end{bmatrix}
   =
   P_G
   \begin{bmatrix}
   u_x \\ u_w
   \end{bmatrix},
   \qquad
   P_G = Z  \widehat{P}^{-1} Z^T,
   \quad
   \begin{bmatrix}
    \bar{u}_x \\ \bar{u}_w
   \end{bmatrix}
   = Z \, \widecheck u ,
   \quad
   Z^T
   \begin{bmatrix}
   u_x \\ u_w
   \end{bmatrix}
   = \widehat u.
\end{equation}
% where
%\begin{equation}
%\label{eq:proj_obl}
% $$
%       P_G =
%       Z  \widehat{P}^{-1} Z^T,
% $$
%       =
%        \begin{bmatrix}
%         G & \\
%             & F^{-1}
%       \end{bmatrix},
%\end{equation}
%
% $$
%     %\label{eq:proj-oper}
%     P_G = Z  \widehat{P}^{-1} Z^T
%     %\begin{bmatrix}
%     %  Z_1 \\ Z_2
%     %\end{bmatrix}
%     %\left(
%       %\begin{bmatrix}
%       %  Z_1^T & Z_2^T
%       %\end{bmatrix}
%       %\begin{bmatrix}
%       % G & \\
%       %    & F^{-1}
%       %\end{bmatrix}
%       %\begin{bmatrix}
%       %  Z_1 \\ Z_2
%       %\end{bmatrix}
%     %\right)^{-1}
%     %\begin{bmatrix}
%     %  Z_1^T & Z_2^T
%     %\end{bmatrix} ,
% $$
% and
% $$
%  \begin{bmatrix}
%       v_x \\ v_w
%  \end{bmatrix}
%  = Z \, \widehat v ,
%  \quad
%  Z^T
%  \begin{bmatrix}
%  u_x \\ u_w
%  \end{bmatrix}
%  = \widehat u .
% $$
Furthermore, %it is easy to verify that %~\eqref{eq:proj_obl}  % nothing here is easy
$$
    P_G
    \begin{bmatrix} G & \\ & F^{-1} \end{bmatrix}
$$
% \smarttodo{{\color{red} The previous matrix, and not $P_G^{-1}$, is a projector.}}
is an oblique projector into $\Null(N)$.
Let $\widehat{L}$ be the lower triangular Cholesky factor of $\widehat{P}$ and let
\begin{equation}
\label{eq:red-kryl-space}
   \widehat{\mathcal{K}} = \mathcal{K} \left( \widehat{L}^{-1} \widehat{M} \, \widehat{L}^{-T}, \;
             \widehat{L}^{-1} (\widehat{b} - \widehat{M} \,\widehat{x}_0) \right)
\end{equation}
be the Krylov space generated by the preconditioned reduced operator $\widehat{L}^{-1} \widehat{M} \, \widehat{L}^{-T}$ and initial vector $\widehat{L}^{-1} (\widehat{b} - \widehat{M} \,\widehat{x}_0)$,
where $\widehat{b}$ is given in~\eqref{eq:reduced-sys-2} and $x_0 = Z_1 \widehat{x}_0$, with $Z_1$ defined in~\eqref{eq:nullspace-basis}.

The computation of~\eqref{eq:apply-PG} can be obtained by solving
\begin{equation}
   \label{eq:proj}
   \begin{bmatrix}
      G &                  & B^T \\
        & F^{-1}           & E^T \\
      B & E\phantom{^{-1}} &
   \end{bmatrix}
   \begin{bmatrix}
      \bar{u}_x \\ \bar{u}_w \\ \bar{z}
   \end{bmatrix}
   =
   \begin{bmatrix}
      u_x \\ u_w \\ 0
   \end{bmatrix},
\end{equation}
see, e.g., \cite{gould-hribar-nocedal-2001}, so that \(P_G\) could be expressed as
\[
  P_G =
  \begin{bmatrix}
       I & 0 & 0
    \\ 0 & I & 0
  \end{bmatrix}
  \begin{bmatrix}
     G &                  & B^T \\
       & F^{-1}           & E^T \\
     B & E\phantom{^{-1}} &
  \end{bmatrix}^{-1}
  \begin{bmatrix}
       I & 0
    \\ 0 & I
    \\ 0 & 0
  \end{bmatrix}.
\]

We now apply Principles~2.1 and~2.2 of \cite{gould-orban-rees-2014} to the standard Lanczos basis-generation
process for $\widehat{\mathcal{K}}$, and obtain the projected Lanczos process outlined in \Cref{alg:proj-lanczos}.
% Although the Lanczos process can be considered as a special case of the Arnoldi process, we report both for the sake of clarity and completeness.
% Note that we highlight the components of the vectors $u_k$ and $v_k$ corresponding to the first two blocks of rows and columns of~\eqref{eq:rsp-2},
% which will be used later.

In \Cref{alg:proj-lanczos}, the notation \(\|u\|_{[P]}\) represents a measure of the deviation of \(u = [u_x \, ; \, u_w]\) from \(\Null(N)\)
\citep[Section~3]{gould-orban-rees-2014}.
More precisely
\begin{equation}
  \label{eq:def-Pnorm}
  \|u\|_{[P]}^2 := u_x^T \bar{u}_x + u_w^T \bar{u}_w,
\end{equation}
where \(\bar{u} = [\bar{u}_x \, ; \, \bar{u}_w]\) is defined by~\eqref{eq:proj}.
Note that \(\|u\|_{[P]}\) is actually a seminorm and vanishes if and only if \([u_x \,;\, u_w]\) is orthogonal to \(\Null(N)\).

\begin{algorithm}[ht]
    \caption{Projected Lanczos Process \label{alg:proj-lanczos}} % for~\eqref{eq:rsp-2}}
    \begin{algorithmic}[1]
      \State choose $[x_0 \,;\, w_0]$ such that \(B x_0 + E w_0 = 0\) \Comment{initial guess}
      \State \label{alg:l-line-v0}%
             $v_{0,x} = 0$, \ $v_{0,w} = -w_0$
             \Comment{initial Lanczos vector}
      \State \label{alg:l-line-proj-u0}%
             $u_{0,x} = b - A x_0$, \ $u_{0,w} = -F^{-1} w_0$
             \Comment{$u_0 = b_0 - M g_0$}
      \State \label{alg:l-line-prec-u0}%
             $[\bar{u}_{1,x} \,;\, \bar{u}_{1,w} \,;\, \bar{z}_1] \leftarrow$ solution of~\eqref{eq:proj} with right-hand side $[u_{0,x} \,;\, u_{0,w} \,;\, 0]$
      \State \label{alg:l-line-proj-v1}%
             $v_{1,x} = \bar{u}_{1,x}$, \ $v_{1,w} = \bar{u}_{1,w}$
             \Comment{$v_1 = P_G \, u_0$}
      \State \label{alg:l-line-proj-beta1}%
             $\beta_1 = (v_{1,x}^T u_{0,x} + v_{1,w}^T u_{0,w})^{\tfrac{1}{2}}$
             \Comment{$\beta_1 = (v_1^T u_0)^{\tfrac{1}{2}}$}
%      \IfOneLine {$\beta_1 \neq 0$}  \label{alg:l-line-scale-v1} \ $v_{1,x} = v_{1,x} / \beta_{1}$, \ $v_{1,w} = v_{1,w} / \beta_{1}$
      \If {$\beta_1 \neq 0$}
             \State \label{alg:l-line-scale-v1}
             $v_{1,x} = v_{1,x} / \beta_{1}$, \ $v_{1,w} = v_{1,w} / \beta_{1}$
             \Comment{$\|v_1\|_{[P]} = 1$}
      \EndIf
      \State $k = 1$
      \While {$\beta_k \ne 0$}
          \State \label{alg:l-line-uk}%
                 $u_{k,x} = A v_{k,x}$, \ $u_{k,w} = F^{-1} v_{k,w}$
                 \Comment{$u_k = M v_k$}
          \State \label{alg:l-line-alphak}%
                 $\alpha_k = v_{k,x}^T u_{k,x} + v_{k,w}^T u_{k,w}$
                 \Comment{$\alpha_k = v_k^T u_k$}
          \State \label{alg:l-line-prec-uk}%
               $[\bar{u}_{k+1,x} \,;\, \bar{u}_{k+1,w} \,;\, \bar{z}_{k+1}] \leftarrow$ solution of~\eqref{eq:proj} with right-hand side $[u_{k,x} \,;\, u_{k,w} \,;\, 0]$
          \State \label{alg:l-line-proj-vk+1-x}% $v_{k+1} = [v_{k+1,x} \,;\, v_{k+1,w}] = P_G \, u_k - \alpha_k v_k - \beta_k v_{k-1}$
                $v_{k+1,x} = \bar{u}_{k+1,x} - \alpha_k v_{k,x} - \beta_k v_{k-1,x}$
                \Comment{$v_{k+1} = \bar{u}_{k+1} - \alpha_k v_k - \beta_k v_{k-1}$}
          \State \label{alg:l-line-proj-vk+1-w}%
                $v_{k+1,w} = \bar{u}_{k+1,w} - \alpha_k v_{k,w} - \beta_k v_{k-1,w}$
          \State \label{alg:l-line-betak+1}%
                 $\beta_{k+1} = (v_{k+1,x}^T u_{k,x} + v_{k+1,w}^T u_{k,w})^{\tfrac{1}{2}}$
                 \Comment{$\beta_{k+1} = (v_{k+1}^T u_k)^{\tfrac{1}{2}}$}
%          \IfOneLine {$\beta_{k+1} \neq 0$} \label{alg:l-line-scale-vk+1} \ $v_{k+1,x} = v_{k+1,x} / \beta_{k+1}$, \ $v_{k+1,w} = v_{k+1,w} / \beta_{k+1}$
          \If {$\beta_{k+1} \neq 0$}
                 \State \label{alg:l-line-scale-vk+1}
                 $v_{k+1,x} = v_{k+1,x} / \beta_{k+1}$, \ $v_{k+1,w} = v_{k+1,w} / \beta_{k+1}$
                 \Comment{$\|v_{k+1}\|_{[P]} = 1$}
          \EndIf
         \State $k = k+1$
      \EndWhile
    \end{algorithmic}
 \end{algorithm}

%\smarttodo{In the comments of \Cref{alg:proj-lanczos} you wrote $\|v_1\|_{\widehat{P}^{-1}} = 1$ and $\|v_{k+1}\|_{\widehat{P}^{-1}} = 1$.
%Is it correct? Check, e.g., the dimensions of $\widehat{P}$ and $v_{k+1}$.}
%\smarttodo[color=white,linecolor=gray]{Sorry I read too fast. It should be good now. It's the same [P]-norm that we use on p.11}

Conceptually, the Lanczos process corresponding to \Cref{alg:proj-lanczos} can be summarized as
\begin{equation*}
  % \label{eq:lanczos-3x3}
  \begin{bmatrix}
     A &                  & B^T \\
       & F^{-1}           & E^T \\
     B & E\phantom{^{-1}} &
  \end{bmatrix}
  \begin{bmatrix}
       V_{k,x}
    \\ V_{k,w}
    \\ \bar{Z}_k
  \end{bmatrix}
  =
  \begin{bmatrix}
     G &                  & B^T \\
       & F^{-1}           & E^T \\
     B & E\phantom{^{-1}} &
  \end{bmatrix}
  \left(
  \begin{bmatrix}
       V_{k,x}
    \\ V_{k,w}
    \\ \bar{Z}_k
  \end{bmatrix}
  T_k
  + \beta_{k+1}
  % \begin{bmatrix}
  %    G &                  & B^T \\
  %      & F^{-1}           & E^T \\
  %    B & E\phantom{^{-1}} &
  % \end{bmatrix}
  \begin{bmatrix}
       v_{k+1,x}
    \\ v_{k+1,w}
    \\ \bar{z}_{k+1}
  \end{bmatrix}
  e_k^T
  \right),
\end{equation*}
where
\[
  V_{k,x} =
  \begin{bmatrix}
    v_{1,x} & \dots & v_{k,x}
  \end{bmatrix},
  \quad
  V_{k,w} =
  \begin{bmatrix}
    v_{1,w} & \dots & v_{k,w}
  \end{bmatrix},
  \quad
  \bar{Z}_k =
  \begin{bmatrix}
    \bar{z}_1 & \dots & \bar{z}_k
  \end{bmatrix},
\]
and \(T_k\) is the usual Lanczos tridiagonal matrix.
Provided that \([x_0 \,;\, w_0] \in \Null(N)\), \cite[Theorem~$2.2$]{gould-orban-rees-2014} guarantees that \Cref{alg:proj-lanczos} is well defined and equivalent
to \Cref{alg:prec-lanczos} in~\Cref{sec:processes} applied to~\eqref{eq:reduced-sys}--\eqref{eq:reduced-sys-2} with preconditioner~\eqref{eq:reduced-prec}.
In \Cref{alg:proj-lanczos} and subsequent algorithms, we use the symbol ``\(\leftarrow\)'' to assign to the vector on the left of the arrow the result of the external
procedure on the right of the arrow.
%, or to indicate that a vector or scalar on the left is overwritten with the value on the right.

In the next sections we show how the projected basis-generation procedures can be further
reformulated by referring to the original system~\eqref{eq:rsp}, thus avoiding the use of $E$ and $F$ and the factorization~\eqref{eq:factc}.
% This leads to the design of constraint-preconditioned variants of these procedures and of the associated Krylov solvers, which is the objective of this work.

\section{Constraint-Preconditioned Lanczos Process\label{sec:cp-lanczos}}

% In \Cref{alg:proj-lanczos}, the computation of $v_k = [v_{k,x} \,;\, v_{k,w}]$ at lines~\ref{alg:l-line-proj-v1} and~\ref{alg:l-line-proj-vk+1} requires the computation of $\bar{u}_{k+1} = P_G u_k$, which can be obtained by solving
% \begin{equation}
%    \label{eq:proj}
%    \begin{bmatrix}
%       G &                  & B^T \\
%         & F^{-1}           & E^T \\
%       B & E\phantom{^{-1}} &
%    \end{bmatrix}
%    \begin{bmatrix}
%       \bar{u}_{k+1,x} \\ \bar{u}_{k+1,w} \\ \bar{z}_{k+1}
%    \end{bmatrix}
%    =
%    \begin{bmatrix}
%       u_{k,x} \\ u_{k,w} \\ 0
%    \end{bmatrix},
% \end{equation}
% see, e.g., \cite{gould-hribar-nocedal-2001}, so that \(P_G\) could be expressed as
% \[
%   P_G =
%   \begin{bmatrix}
%        I & 0 & 0
%     \\ 0 & I & 0
%   \end{bmatrix}
%   \begin{bmatrix}
%      G &                  & B^T \\
%        & F^{-1}           & E^T \\
%      B & E\phantom{^{-1}} &
%   \end{bmatrix}^{-1}
%   \begin{bmatrix}
%        I & 0
%     \\ 0 & I
%     \\ 0 & 0
%   \end{bmatrix}.
% \]
If we define \(\bar{p}_k = \bar{u}_{k,x}\) for all \(k \geq 1\), and
\begin{equation}
\label{eq:tk}
  t_k = EFu_{k,w}, \quad k = 0, 1, \dots
\end{equation}
then~\eqref{eq:proj} at line~\ref{alg:l-line-prec-uk} of \Cref{alg:proj-lanczos} can be written as
\begin{equation}
  \label{eq:proj-2}
  \begin{bmatrix}
   G & \phantom{-}B^T \\
   B & -C\phantom{^T}
  \end{bmatrix}
  \begin{bmatrix}
    \bar{p}_{k+1} \\ \bar{z}_{k+1}
  \end{bmatrix}
  =
  \begin{bmatrix}
     u_{k,x} \\ -t_k
  \end{bmatrix}.
\end{equation}
\Cref{assum:Phat-spd} occurs when the sum of the number of negative eigenvalues of the matrix of~\eqref{eq:proj-2} and \(C\) is \(m\) \cite[Theorem~$2.1$]{dollar-gould-schilders-wathen-2006}, which may be verified if an inertia-revealing symmetric indefinite factorization is used to solve~\eqref{eq:proj-2}, such as that of \cite{duff-2004}.

Unfortunately, \eqref{eq:proj-2} still appears to depend on $F$ via~\eqref{eq:tk}.
We now reformulate \Cref{alg:proj-lanczos} in terms of full-space quantities.
% In the discussion below, all line numbers refer to \Cref{alg:proj-lanczos}.
Define the initial guess
\begin{equation}
\label{eq:w0q0}
   w_0 = -FE^T q_0,
\end{equation}
where $q_0 \in \R^m$ is arbitrary (e.g., $q_0=0$).
Line~\ref{alg:l-line-proj-u0} of \Cref{alg:proj-lanczos} and~\eqref{eq:tk}  yield
\[
   u_{0,x} = r_0 = b - A x_0, \quad
   u_{0,w} = E^T q_0, \quad
   t_0 = C q_0.
\]
From here on, let us denote \(p_k = v_{k,x}\).
At lines~\ref{alg:l-line-prec-u0}-\ref{alg:l-line-proj-v1} of \Cref{alg:proj-lanczos}, we compute \(p_{1} = \bar{p}_{1}\) and \(\bar{z}_1\) from~\eqref{eq:proj}, which yields, in particular, \(\bar{u}_{1,w} = F E^T (q_0 - \bar{z}_1)\).
If we define
\[
   s_1 = q_0 - \bar{z}_1, \quad
   q_1 = s_1,
\]
lines~\ref{alg:l-line-proj-v1}-\ref{alg:l-line-proj-beta1} of \Cref{alg:proj-lanczos} take the form
\begin{subequations}
  \label{eq:w1q1}
  \begin{align}
    p_{1}   & = \bar{p}_{1} \label{eq:w1q1-v1x} \\
    v_{1,w} & = FE^T q_1, \label{eq:w1q1-v1w} \\
    \beta_1 & = ( p_{1}^T u_{0,x} + q_1^T C q_1 )^{\tfrac{1}{2}} = ( p_{1}^T u_{0,x} + q_1^T t_0 )^{\tfrac{1}{2}}. \label{eq:w1q1-beta1}
  \end{align}
\end{subequations}
% where we used the second block of~\eqref{eq:proj} to obtain $v_{1,w}$.
We then normalize by dividing $p_{1}$ and $q_1$ by $\beta_1$.
Lines~\ref{alg:l-line-uk}--\ref{alg:l-line-alphak} of \Cref{alg:proj-lanczos} and~\eqref{eq:tk} give
\begin{subequations}
  \begin{align}
    u_{1,w}  & = E^T q_1, \label{eq:u1w}  \\
    t_1      & =  C q_1,  \label{eq:t1} \\
    \alpha_1 & = p_{1}^T u_{1,x} + q_1^T t_1 = p_{1}^T A p_{1} + q_1^T C q_1. \label{eq:alpha1}
  \end{align}
\end{subequations}
We now compute \(p_2\) and \(v_{2,w}\) from lines~\ref{alg:l-line-proj-vk+1-x}--\ref{alg:l-line-proj-vk+1-w} of \Cref{alg:proj-lanczos} with $k=1$, i.e., we compute \(\bar{p}_{2}\) and \(\bar{z}_2\) from~\eqref{eq:proj-2}, and note that \(\bar{u}_{2,w} = F E^T (q_1 - \bar{z}_2)\).
Thus, by setting
\[
  s_2 = q_1 - \bar{z}_2, \quad
  q_2 = s_2 - \alpha_1 q_1 - \beta_1 q_0,
\]
we obtain from lines~\ref{alg:l-line-v0} and \ref{alg:l-line-proj-vk+1-x}--\ref{alg:l-line-proj-vk+1-w} of \Cref{alg:proj-lanczos} together with~\eqref{eq:w1q1-v1w}, \eqref{eq:u1w} and~\eqref{eq:t1}:
\begin{align*}
   p_{2}   & = \bar{p}_{2} - \alpha_1 p_{1} - \beta_1 p_{0} \\
   v_{2,w} & = FE^T s_2 - \alpha_1 FE^T q_1 - \beta_1 FE^T q_0 = FE^T q_2 \\
   \beta_2 & = (p_{2}^T A p_{1} + q_2^T C q_1)^{\tfrac{1}{2}} = (p_{2}^T u_{1,x} + q_2^T t_1 )^{\tfrac{1}{2}}.
\end{align*}
% where $v_{2,x}$ and $v_{2,w}$ appear in line~\ref{alg:l-line-proj-vk+1} of \Cref{alg:proj-lanczos} for $k=1$.
Then, according to line~\ref{alg:l-line-scale-vk+1},
$p_{2}$ must be divided by $\beta_2$, and we do the same with $q_2$.
An induction argument shows that for all $k \ge 1$
\begin{align*}
   u_{k,w} & = E^T q_k,  \\
   t_k & = C q_k, \\
   \alpha_k & = p_{k}^T u_{k,x} + q_k^T t_k = p_{k}^T A p_{k} + q_k^T C q_k,
\end{align*}
where $q_k$ has been normalized by $\beta_k$.
Furthermore, letting
\[
   s_{k+1} = q_k - \bar{z}_{k+1}, \quad
   q_{k+1} = s_{k+1} - \alpha_k q_k - \beta_k q_{k-1},
\]
we obtain
\begin{align*}
    p_{k+1}     & = \bar{p}_{k+1} - \alpha_k p_{k} - \beta_k p_{k-1} \\
    v_{k+1,w}   & = FE^T q_{k+1}, \\
    \beta_{k+1} & = (p_{k+1}^T A p_{k} + q_{k+1}^T C q_k)^{\tfrac{1}{2}} = (p_{k+1}^T u_{k,x} + q_{k+1}^T t_k )^{\tfrac{1}{2}}.
\end{align*}
We divide $p_{k+1}$ and $q_{k+1}$ by $\beta_{k+1}$ to obtain the vectors to be used at the next iteration.
Thus, if we rename $u_{k,x}$ as $u_k$, we obtain \Cref{alg:cp-lanczos}.

\begin{algorithm}[ht]
    \caption{Constraint-Preconditioned Lanczos Process \label{alg:cp-lanczos}}
                                 % for \eqref{eq:rsp}}
    \begin{algorithmic}[1]
      \State choose $[x_0 \,;\, q_0]$ such that \(B x_0 - C q_0 = 0\) \Comment{initial guess}
      \State $p_0 = 0$ \Comment{initial Lanczos vector}
      \State $u_0 = b - A x_0$, \ $t_0 = C q_0$       \label{alg:cpl-line-u0}
      \State \label{alg:cpl-line-proj2-0}%
             $[\bar{p}_1 \,;\, \bar{z}_1] \leftarrow$ solution of~\eqref{eq:proj-2} with right-hand side $[u_0 \,;\, -t_0]$
      \State $p_1 = \bar{p}_1$                          \label{alg:cpl-line-p1}
      \State $s_1 = q_0 - \bar{z}_1$, \ $q_1 = s_1$    \label{alg:cpl-line-q1}
      \State $\beta_1 = (p_1^T u_0 + q_1^T t_0)^{\tfrac{1}{2}}$   \label{alg:cpl-line-beta1}
%      \IfOneLine {$\beta_1 \ne 0$} \label{alg:cpl-line-scale1} \ $p_1 = p_1 / \beta_1$, \ $q_1 = q_1 / \beta_1$
      \If {$\beta_1 \ne 0$}
            \State \label{alg:cpl-line-scale1} $p_1 = p_1 / \beta_1$, \ $q_1 = q_1 / \beta_1$
      \EndIf
      \State $k = 1$
      \While {$\beta_k \ne 0$}
         \State $u_k = A p_k$, \ $t_k = C q_k$
         \State $\alpha_k = p_k^T u_k + q_k^T t_k$ \Comment{$= p_k^T A p_k + q_k^T C q_k$}   \label{alg:cpl-line-alphak}
         \State \label{alg:cpl-line-proj2-k}%
                $[\bar{p}_{k+1} \,;\, \bar{z}_{k+1}] \leftarrow$ solution of~\eqref{eq:proj-2} with right-hand side $[u_k \,;\, -t_k]$
         \State $p_{k+1} = \bar{p}_{k+1} - \alpha_k p_k - \beta_k p_{k-1}$ \label{alg:cpl-line-pk+1}
         \State $s_{k+1} = q_k - \bar{z}_{k+1}$, \ $q_{k+1} = s_{k+1} - \alpha_k q_k - \beta_k q_{k-1}$
                                                         \label{alg:cpl-line-qk+1}
          \State $\beta_{k+1} = (p_{k+1}^T u_k + q_{k+1}^T t_k)^{\tfrac{1}{2}}$ \Comment{$= (p_{k+1}^T A p_k + q_{k+1}^T C q_k)^{\tfrac{1}{2}}$}
%        \IfOneLine {$\beta_{k+1} \ne 0$} \label{alg:cpl-line-scalek+1} \ $p_{k+1} = p_{k+1} / \beta_{k+1}$, \ $q_{k+1} = q_{k+1} / \beta_{k+1}$
          \If {$\beta_{k+1} \ne 0$}
                \State \label{alg:cpl-line-scalek+1} $p_{k+1} = p_{k+1} / \beta_{k+1}$, \ $q_{k+1} = q_{k+1} / \beta_{k+1}$
          \EndIf
          \State $k = k+1$
       \EndWhile
    \end{algorithmic}
\end{algorithm}

The above transformations can be condensed in the following principle, which summarizes the conversion a of projected process into a constraint-preconditioned process.

\begin{principle}
  \label{pri:proj-to-cp}
  \mbox{~}
  \begin{enumerate}
    \item Basis vectors \(v_{k+1,x}\) are unchanged;
    \item Basis vectors \(v_{k+1,w}\) have the form \(FE^T q_{k+1}\), where \(q_{k+1}\) is defined by
    \begin{align*}
      s_{k+1} & = q_k - \bar{z}_{k+1}, \\
      q_1       & = s_1, \\
      q_{k+1} & = s_{k+1} - \alpha_k q_k - \beta_k q_{k-1}, \quad (k \geq 1),
    \end{align*}
    and where \(\bar{z}_{k+1}\) results from the solution of~\eqref{eq:proj-2};
    \item Inner products of the form \(v_{i,w}^T u_{j,w}\) become \(q_i^T C q_j = q_i^T t_j\).
  \end{enumerate}
\end{principle}

\Cref{thm:cp-lanczos-equiv} summarizes the equivalence between the two formulations.

\begin{shadytheorem}
	\label{thm:cp-lanczos-equiv}
  Let $E$ and $F$ be as defined in~\eqref{eq:factc} and $G$ chosen to satisfy \Cref{assum:Phat-spd}.
  Let $q_0 \in \R^m$ be arbitrary. Then, \Cref{alg:proj-lanczos} with starting guesses $x_0 \in \R^n$ and $w_0 = -FE^T q_0$,
  such that $B x_0 + E w_0$, is equivalent to \Cref{alg:cp-lanczos}
  with starting guesses $x_0$ and $q_0$. In particular, for all~$k$, the vectors $v_{k,x}$ and $v_{k,w}$, and the scalars $\alpha_k$
  and $\beta_k$ in \Cref{alg:proj-lanczos} are equal to the vectors $p_k$ and $FE^T q_k$, and to the scalars $\alpha_k$ and
  $\beta_k$ in \Cref{alg:cp-lanczos}, respectively.
\end{shadytheorem}

\smallskip

Note that \Cref{alg:cp-lanczos} does not contain references to $E$ and $F$.
The variable $s_k$ is used only to improve readability.
\Cref{assum:Phat-spd} guarantees that \Cref{alg:cp-lanczos} is well posed because it is equivalent to \Cref{alg:proj-lanczos}, which, in turn, is equivalent to the standard Lanczos process for building an orthonormal basis of~\eqref{eq:red-kryl-space}.
The main advantages of \Cref{alg:cp-lanczos} are that it works directly with the formulation~\eqref{eq:rsp} and it only requires storage for three vectors of size \(n + m\) ($[p_k \, ; \, q_k]$, $[u_k \, ; \, t_k]$, and $[\bar{p}_k \, ; \, \bar{z}_k]$), as opposed to the same number of vectors of size \(n + p + m\) for \Cref{alg:proj-lanczos}.

We call \Cref{alg:cp-lanczos} the Constraint-Preconditioned Lanczos (CP-Lanczos) process because of its similarity to a Lanczos process for building an orthonormal basis of a Krylov space associated with the preconditioned operator $P^{-1} M$, even though the latter appears nonsymmetric.

% In the next section we show how the CP-Lanczos process can be used to obtain Krylov methods for regularized saddle-point systems with symmetric leading block.
% As we will see, the choice of the starting guess cannot be completely arbitrary.

%the values of $\alpha_k$ and $\beta_k$ computed by this algorithm are the same
%as the ones computed by \Cref{alg:proj-lanczos}, and, according to the
%principles set out in~\cite{gould-orban-rees-2014}, are the entries of the
%matrix
%$$
%    T_{k,k+1} =
%    \begin{bmatrix}
%    \alpha_1 & \beta_1  &             &              &             \\
%    \beta_1  & \alpha_2 & \beta_2     &              &             \\
%             & \ddots   & \ddots      & \ddots       &             \\
%             &          & \beta_{k-2} & \alpha_{k-1} & \beta_{k-1} \\
%             &          &             & \beta_{k-1}  & \alpha_k    \\
%             &          &             &              & \beta_k
%    \end{bmatrix}
%$$

\section{Constraint-Preconditioned Lanczos-Based Krylov Solvers\label{sec:cp-l-methods}}

We may exploit \Cref{thm:cp-lanczos-equiv} and use \Cref{alg:cp-lanczos} to derive a
constraint-preconditioned version of any Krylov method based on the Lanczos process.
To this aim, we must understand how the update of the $k$-th iterate $[x_k \,;\, w_k]$ in a Krylov method based on \Cref{alg:proj-lanczos}
translates into the update of the $k$-th iterate $[x_k \,;\, y_k]$ in the version of that Krylov method based on \Cref{alg:cp-lanczos}.
In the following, the former and the latter version of the Krylov method are referred to as projected-Krylov (P-Krylov) and
constraint-preconditioned-Krylov (CP-Krylov), respectively.
% Of course, we must choose a suitable starting guess for the CP-Krylov method.

Because the initial guess $g_0 = [x_0 \,;\, w_0]$ of P-Krylov applied to~\eqref{eq:rsp-2}
must lie in $\Null(N)$, CP-Krylov must be initialized with $[x_0 \,;\, y_0]$ such that
\begin{equation}
\label{eq:2ndeq0_x0y0}
   B x_0 - C y_0 = 0.
\end{equation}
% Furthermore, from~\eqref{eq:w} and~\eqref{eq:w0q0} it follows that we can set $q_0 = -y_0$.
% \smarttodo{I find this part, from here until the theorem, very confusing and hard to follow. I am not sure it is necessary.}
% \smarttodo{{\color{red} Check if we must set $q_0 = y_0 = 0$ and $x_0 = 0$. Note that $q_0$ is
% both a fake Lanczos vector and starting guess.}}
% Because $[\bar{p}_{k+1}  \,;\, \bar{z}_{k+1}]$ solves~\eqref{eq:proj-2}, the second block equation yields
% \begin{equation}
%    \label{eq:2ndeq0-barpkqk}
%    B \bar{p}_{k+1} - C (\bar{z}_{k+1} - q_k) = 0.
% \end{equation}
% By comparing line~\ref{alg:l-line-proj-vk+1} of \Cref{alg:proj-lanczos} with
% line~\ref{alg:cpl-line-qk+1} of \Cref{alg:cp-lanczos}, we see that the
% $(k+1)$-st Lanczos vector computed by CP-Krylov is
% \smarttodo{Is it $-q_{k+1}$ or $q_{k+1}$? {\color{red} It should be $-q_{k+1}$.} I don't understand this.}
% $[p_{k+1} \,;\, -q_{k+1}]$.
% Thus, the update of $[x_{k+1} \,;\, y_{k+1}]$ in CP-Krylov has the same
% form as the update of $[x_{k+1} \,;\, w_{k+1}]$ in P-Krylov, provided that
% $[p_{k+1} \,;\, -q_{k+1}]$ is used instead of $[v_{k+1,x} \,;\, v_{k+1,w}]$.
% Furthermore, $[x_{k+1} \,;\, w_{k+1}] \in  \Null([B \;\, E])$ implies
% \begin{equation}
% \label{eq:2ndeq0_xkyk}
%    B x_{k+1} - C y_{k+1} = 0.
% \end{equation}
% It is also easy to show that
% \smarttodo{$-C q_{k+1}$? {\color{red} It should be $C q_{k+1}$.} With $k=0$, we have $B p_1 + C q_1 = t_0 \neq 0$.}
% \begin{equation}
% \label{eq:2ndeq0_pkqk}
%    B p_{k+1} + C q_{k+1} = 0
% \end{equation}
% and
% $$
%     y_{k+1} = y_0 - q_{k+1}.
% $$

Our first result states a property of \Cref{alg:cp-lanczos} that follows from a specific \(q_0\).

\begin{shadylemma}
  \label{lem:cp-lanczos}
  Let \Cref{alg:cp-lanczos} be initialized with \(x_0 \in \R^n\) and \(q_0 \in \Null(C)\).
  Then, for all \(k \geq 0\),
  \begin{equation}
    \label{eq:cp-lanczos-Bp+Cq=0}
    B p_k + C q_k = 0.
  \end{equation}
\end{shadylemma}

\begin{proof}
  We proceed by induction.
  For \(k = 0\),~\eqref{eq:cp-lanczos-Bp+Cq=0} holds because \(p_0 = 0\) and \(q_0 \in \Null(C)\).
  For \(k = 1\), \(p_1 = \bar{p}_1\), \(q_1 = q_0 - \bar{z}_1 = -\bar{z}_1\), and~\eqref{eq:proj-2} and our assumption that \(q_0 \in \Null(C)\) yield
  \[
    B p_1 + C q_1 = B \bar{p}_1 - C \bar{z}_1 = -t_0 = -C q_0 = 0.
  \]
  Assume~\eqref{eq:cp-lanczos-Bp+Cq=0} holds for any index \(j \le k\).
  Lines~\ref{alg:cpl-line-pk+1}--\ref{alg:cpl-line-qk+1} of \Cref{alg:cp-lanczos}, \eqref{eq:tk}, \eqref{eq:proj-2}, and our induction assumption imply that
  \begin{align*}
       B p_{k+1} + C q_{k+1}
       & = B \bar{p}_{k+1} + C (q_k - \bar{z}_{k+1}) - \alpha_k (B p_k + C q_k) - \beta_k (B p_{k-1} + C q_{k-1})
    \\ & = B \bar{p}_{k+1} + C (q_k - \bar{z}_{k+1})
    \\ & = B \bar{p}_{k+1} - C \bar{z}_{k+1} + t_k = 0,
  \end{align*}
  which establishes~\eqref{eq:cp-lanczos-Bp+Cq=0}.
\end{proof}

An interesting property of the CP-Lanczos process is that it is equivalent to formally applying the standard Lanczos process
to system~\eqref{eq:rsp} with preconditioner~\eqref{eq:cp}, where by ``formal application'', we mean that the Lanczos process is applied blindly
as if $P$ were positive definite. Such formal application is stated as \Cref{alg:cp-lanczos-full-space} in \Cref{sec:processes}.
The equivalence with \Cref{alg:cp-lanczos} is stated in the next result, which parallels \cite[Theorem~$2.2$]{gould-orban-rees-2014}.

\begin{shadytheorem}
  \label{thm:cp-lanczos-equiv-2}
  Let \Cref{alg:cp-lanczos} be initialized with \(x_0 \in \R^n\) such that \(B x_0 = 0\), \(q_0 = 0 \in \R^m \),
  and \Cref{alg:cp-lanczos-full-space} be initialized with the same \(x_0\) and \(y_0 \in \R^m\) such
  that~\eqref{eq:2ndeq0_x0y0} is satisfied. Then, for all \(k \geq 0\), \(v_{k,x} = p_k\) and \(v_{k,y} = -q_k\),
  where \([v_{k,x} \,;\, v_{k,y}]\) is the \(k\)-th Lanczos vector generated in \Cref{alg:cp-lanczos-full-space},
  and \(p_k\) and \(q_k\) are the \(k\)-th Lanczos vectors generated in \Cref{alg:cp-lanczos}.
  In addition, the scalars \(\alpha_k\) and \(\beta_k\) computed at each iteration are the same in both algorithms.
\end{shadytheorem}

\begin{proof}
  We proceed by induction.
  The result holds for \(k = 0\) because \([v_{0,x} \,;\, v_{0,y}] = [0 \,;\, 0] = [p_0 \,;\, -q_0]\).
  With \(q_0 = 0\), \Cref{alg:cp-lanczos} initializes \(u_0 = b - A x_0\) and \(t_0 = 0\).
  Because~\eqref{eq:2ndeq0_x0y0} is satisfied, \Cref{alg:cp-lanczos-full-space} initializes \(r_{0,x} = u_0 - B^T y_0\) and \(r_{0,y} = 0\).
  Thus, \([v_{1,x} \,;\, v_{1,y}]\) solves~\eqref{eq:proj-2} with right-hand side \([u_0 - B^T y_0 \,;\, 0]\).
  By \cite[Theorem~$2.1$, item~$2$]{gould-orban-rees-2014}, \([v_{1,x} \,;\, v_{1,y}]\) equivalently solves~\eqref{eq:proj-2} with right-hand
  side \([u_0 \,;\, 0]\), and therefore, \([v_{1,x} \,;\, v_{1,y}]\) at line~\ref{line:v1-full-space} of \Cref{alg:cp-lanczos-full-space} is equal to \([\bar{p}_1 \,;\, \bar{z}_1]\).
  Lines~\ref{alg:cpl-line-p1}--\ref{alg:cpl-line-q1} of \Cref{alg:cp-lanczos} subsequently set \(p_1 = \bar{p}_1 = v_{1,x}\)
  and \(q_1 = s_1 = q_0 - \bar{z}_1 = -v_{1,y}\).

  With \(q_0 = 0\), line~\ref{alg:cpl-line-beta1} of \Cref{alg:cp-lanczos} computes \(\beta_1 = (p_1^T u_0)^{\tfrac{1}{2}}\).
  We take the inner product of the second row of~\eqref{eq:proj-2} with \(\bar{z}_1 = -q_1\) and note that \(t_0 = 0\),
  and obtain \(\bar{z}_1^T B p_1 = \bar{z}_1^T C \bar{z}_1 = q_1^T C q_1\).
  Similarly, we take the inner product of the first row of~\eqref{eq:proj-2} with \(p_1\) and substitute \(\bar{z}_1^T B p_1\)
  to obtain \(p_1^T u_0 = p_1^T G p_1 + q_1^T C q_1\), so that \(\beta_1\) is the same as that computed at
  line~\ref{alg:cplfs:beta1} of \Cref{alg:cp-lanczos-full-space}.
  We have established that the result also holds for \(k = 1\).

  At a general iteration \(k\), \Cref{alg:cp-lanczos} sets \(u_k = A p_k\), \(t_k = C q_k\) and computes
  \(\alpha_k = p_k^T u_k + q_k^T t_k = p_k^T A p_k + q_k^T C q_k\).
  By \Cref{lem:cp-lanczos}, \(q_k^T B p_k + q_k^T C q_k = 0\), so that \(\alpha_k = p_k^T A p_k - 2 q_k^T B p_k - q_k^T C q_k\).
  Under the recurrence assumption that \(v_{k,x} = p_k\) and \(v_{k,y} = -q_k\), this expression of \(\alpha_k\) is the same as that
  computed at line~\ref{alg:cplfs:alphak} of \Cref{alg:cp-lanczos-full-space}.

  At line~\ref{alg:cpl-line-proj2-k} of \Cref{alg:cp-lanczos}, we compute \([\bar{p}_{k+1} \,;\, \bar{z}_{k+1}]\) from~\eqref{eq:proj-2},
  or, equivalently, as the solution to
  \[
    \begin{bmatrix}
     G & \phantom{-}B^T \\
     B & -C\phantom{^T}
    \end{bmatrix}
    \begin{bmatrix}
      \bar{p}_{k+1} \\ \bar{z}_{k+1} - q_k
    \end{bmatrix}
    =
    \begin{bmatrix}
       A p_k \\ 0
    \end{bmatrix}.
  \]
  In view of \Cref{lem:cp-lanczos}, our recurrence assumption, and \cite[Theorem~$2.1$, item~$2$]{gould-orban-rees-2014},
  line~\ref{alg:cplfs:vk+1} of \Cref{alg:cp-lanczos-full-space} computes \([v_{k+1,x} \,;\, v_{k+1,y}]\) as the solution to the same system as above.
  Therefore, at that point in each algorithm \(v_{k+1,x} = \bar{p}_{k+1}\) and \(v_{k+1,y} = \bar{z}_{k+1} - q_k = -s_{k+1}\).
  The vector updates at lines~\ref{alg:cpl-line-pk+1}--\ref{alg:cpl-line-qk+1} of \Cref{alg:cp-lanczos} together with those at
  line~\ref{alg:cplfs:update-vk+1} of \Cref{alg:cp-lanczos-full-space} show that \(v_{k+1,x} = p_{k+1}\) and \(v_{k+1,y} = -q_{k+1}\).
  Our recurrence assumption and \Cref{lem:cp-lanczos} yield \(B v_{k,x} - C v_{k,y} = 0\) and \(B v_{k+1,x} - C v_{k+1,y} = 0\).
  Finally, \Cref{alg:cp-lanczos-full-space} sets
  \begin{align*}
    \beta_{k+1}^2 & =
    v_{k+1,x}^T u_{k,x} + v_{k+1,y}^T u_{k,y}
    \\ & =
    v_{k+1,x}^T A v_{k,x} + (B v_{k+1,x} - C v_{k+1,y})^T v_{k,y} + v_{k+1,y}^T B v_{k,x}
    \\ & =
    v_{k+1,x}^T A v_{k,x} + v_{k+1,y}^T C v_{k,x},
  \end{align*}
  which is the same value computed in \Cref{alg:cp-lanczos}.
\end{proof}

\Cref{thm:cp-lanczos-equiv-2} shows that \Cref{alg:cp-lanczos} may be summarized as
\[
  \begin{bmatrix}
   A & \phantom{-}B^T \\
   B & -C\phantom{^T}
  \end{bmatrix}
  \begin{bmatrix}
    P_k \\ -Q_k
  \end{bmatrix}
  =
  \begin{bmatrix}
   G & \phantom{-}B^T \\
   B & -C\phantom{^T}
  \end{bmatrix}
  \left(
    \begin{bmatrix}
      P_k \\ -Q_k
    \end{bmatrix}
    T_k
    + \beta_{k+1}
    \begin{bmatrix}
      p_{k+1} \\ -q_{k+1}
    \end{bmatrix}
    e_k^T,
  \right)
\]
provided that \(B x_0 = 0\) and \(q_0 = 0\), where \(T_k\) is the same as in \Cref{alg:proj-lanczos}, and
\[
  P_k =
  \begin{bmatrix}
    p_1 & \dots & p_k
  \end{bmatrix},
  \quad
  Q_k =
  \begin{bmatrix}
    q_1 & \dots & q_k
  \end{bmatrix}.
\]

A consequence of \Cref{thm:cp-lanczos-equiv-2} is that any CP-Krylov method is formally equivalent
to the corresponding standard Krylov method applied to system~\eqref{eq:rsp} with preconditioner~\eqref{eq:cp}.

\begin{shadycorollary}
  \label{cor:p-cp-equivalence}
  Let \Cref{alg:cp-lanczos} be initialized with \(x_0 \in \R^n\) such that \(B x_0 = 0\), \(q_0 = 0 \in \R^m\),
  and \Cref{alg:cp-lanczos-full-space} be initialized with the same \(x_0\) and \(y_0 \in \R^m\) such that~\eqref{eq:2ndeq0_x0y0}
  is satisfied. The $k$-th approximate solution of~\eqref{eq:rsp} computed by any Lanczos-based CP-Krylov method
  coincides with the $k$-th approximate solution obtained by formally applying the standard version of the same method
  to~\eqref{eq:rsp} with preconditioner~\eqref{eq:cp}.
\end{shadycorollary}

Although \Cref{cor:p-cp-equivalence} states that standard Lanczos-based methods can be safely applied to~\eqref{eq:rsp}
with preconditioner~\eqref{eq:cp} and an appropriate starting point, \Cref{alg:cp-lanczos} reduces the computational effort
by never requiring products with \(B\) or \(B^T\). Only products with \(A\) and \(C\) are necessary. On the other hand,
thanks to \Cref{thm:cp-lanczos-equiv-2}, specialized implementations of the standard Lanczos-based methods can be
developed by exploiting the equalities $B p_k + C q_k = 0$ and $B x_k - C y_k = 0$, thus saving matrix-vector products.
% \smarttodo{Daniela: actually, by exploiting the equalities $B p_k + C q_k = 0$ and $B x_k - C y_k = 0$, the standard Lanczos-based methods
% can be implemented by saving many matrix-vector products. Of course, this is a specialized implementation of the standard methods ...}
% \smarttodo[color=white,linecolor=gray]{Dom: That's true. We could mention it in the conclusions, but we don't know if it affects stability}
% \smarttodo{Daniela: from experiments carried out in the past, this improves the stability of standard CG compared with an implementation
% where the previous equalities are not taken into account. However, I added above the observation made in the ``smarttodo'' box.}
% Indeed the matrix-vector product with the whole saddle-point matrix $M$, required at each iteration of standard Krylov methods, is substituted by matrix-vector products involving only the leading and trailing blocks $A$ and $C$;
The computation involving \(s_{k+1}\) can be carried out, for example, as the update \(q_{k-1} = q_k - \bar{z}_{k+1} - \beta_k q_{k-1}\) followed by \(q_{k+1} = q_{k-1} - \alpha_k q_k\), or \(s_{k+1}\) can overwrite \(\bar{z}_{k+1}\).
% furthermore, the additional sum of vectors required to compute $s_k$ in \Cref{alg:cp-lanczos} is generally negligible compared to the cost of performing matrix-vector products with the off-diagonal blocks $B$ and $B^T$ of $M$.
Finally, once~\eqref{eq:cp} has been factorized, storing $B$ is no longer necessary, and this can be used to free memory if needed.

A consequence of \Cref{thm:cp-lanczos-equiv-2} is a formal equivalence between the iterates generated by Lanczos-based
methods applied by way of \Cref{alg:cp-lanczos} and \Cref{alg:cp-lanczos-full-space}.
This equivalence requires a re-interpretation of the optimality conditions associated with
the Krylov method.

Consider, e.g., MINRES~\citep{paige-saunders-1975}.
The residual associated with iterate \([x_k \,;\, w_k \,;\, y_k]\) generated by P-MINRES, with \(w_k = -F E^T y_k\), is
\[
  r_{\text{P},k} =
   \begin{bmatrix}
    r_{\text{P},k,x} \\ r_{\text{P},k,w} \\ r_{\text{P},k,y}
  \end{bmatrix}
  =
  \begin{bmatrix}
    b \\ 0 \\ 0
  \end{bmatrix}
  -
  \begin{bmatrix}
       A &        & B^T
    \\   & F^{-1} & E^T
    \\ B & E      &
  \end{bmatrix}
  \begin{bmatrix}
    x_k \\ w_k \\ y_k
  \end{bmatrix}
  =
  \begin{bmatrix}
       b - A x_k - B^T y_k
    \\ 0
    \\ 0
  \end{bmatrix},
\]
where we used the fact that \(B x_k + E w_k = 0\) for all \(k\).
This residual corresponds to the residual at iterate \([x_k \,;\, y_k]\) generated by CP-MINRES:
\[
  r_{\text{CP},k} =
  \begin{bmatrix}
    b \\ 0
  \end{bmatrix}
  -
  \begin{bmatrix}
    A & \phantom{-}B^T \\
    B & -C
  \end{bmatrix}
  \begin{bmatrix}
    x_k \\ y_k
  \end{bmatrix}
  =
  \begin{bmatrix}
    b - A x_k - B^T y_k
    \\ 0
  \end{bmatrix},
\]
where we exploited the fact that \(B x_k - C y_k = 0\) for all \(k\), which comes from
\(B x_k + E w_k = 0\) and \(w_k = -F E^T y_k\).

%The residual at iterate \((x_k, y_k)\) generated by CP-MINRES is
%\[
%  r_{\text{CP},k} =
%  \begin{bmatrix}
%    b \\ 0
%  \end{bmatrix}
%  -
%  \begin{bmatrix}
%    A & \phantom{-}B^T \\
%    B & -C
%  \end{bmatrix}
%  \begin{bmatrix}
%    x_k \\ y_k
%  \end{bmatrix}
%  =
%  \begin{bmatrix}
%    b - A x_k - B^T y_k
%    \\ 0
%  \end{bmatrix},
%\]
%where we used the fact that \(B x_k - C y_k = 0\) for all \(k\).
%To \((x_k, y_k)\) there corresponds an iterate \((x_k, w_k, y_k)\) generated by P-MINRES such that \(w_k = -F E^T y_k\) and with associated residual
%\[
%  r_{\text{P},k} =
%  \begin{bmatrix}
%    b \\ 0 \\ 0
%  \end{bmatrix}
%  -
%  \begin{bmatrix}
%       A &        & B^T
%    \\   & F^{-1} & E^T
%    \\ B & E      &
%  \end{bmatrix}
%  \begin{bmatrix}
%    x_k \\ w_k \\ y_k
%  \end{bmatrix}
%  =
%  \begin{bmatrix}
%       b - A x_k - B^T y_k
%    \\ 0
%    \\ 0
%  \end{bmatrix},
%\]
%where we used the fact that \(B x_k + E w_k = 0\) for all \(k\).

We may apply the arguments of \citet[Section~3]{gould-orban-rees-2014} to conclude that P-MINRES,
and hence CP-MINRES, minimizes the deviation of \([r_{\text{P},k,x} \,;\, 0]\) from
the range space of \(N\), i.e., as in~\eqref{eq:def-Pnorm},
\begin{equation} \label{eq:norm-dev}
  \|r_{\text{P},k}\|^2_{[P]} = (b - A x_k - B^T y_k)^T h_k,
\end{equation}
where
\[
  \begin{bmatrix}
       G &        & B^T
    \\   & F^{-1} & E^T
    \\ B & E      &
  \end{bmatrix}
  \begin{bmatrix}
       h_k
    \\ f_k
    \\ l_k
  \end{bmatrix}
  =
  \begin{bmatrix}
       b - A x_k - B^T y_k
    \\ 0
    \\ 0
  \end{bmatrix} .
\]
% and, with a little abuse of notation, we use \(r_{\text{P},k}\) instead of \([r_{\text{P},k,x} \,;\, 0]\).
Because \(h_k \in \Null(B)\), we also have
\begin{equation} \label{eq:norm-dev-2}
  \|r_{\text{P},k}\|^2_{[P]} = (b - A x_k)^T h_k.
\end{equation}
Equivalently, \(h_k\) may be computed from
\[
  \begin{bmatrix}
    G & \phantom{-}B^T \\
    B & -C
  \end{bmatrix}
  \begin{bmatrix}
    h_k \\ l_k
  \end{bmatrix}
  =
  \begin{bmatrix}
    b - A x_k
    \\ 0
  \end{bmatrix}.
\]

Because of its residual norm minimization property, CP-MINRES is appropriate to solve saddle-point systems
in a linesearch inexact-Newton context, where we seek to reduce the residual of the Newton-like
equations~\eqref{eq:rsp} in an appropriate space.

The same reasoning applies to the constraint-preconditioned version of any Lanczos-based Krylov method.
For example, \cite{paige-saunders-1975} derive the conjugate gradient method of \cite{hestenes-stiefel-1952} directly from the Lanczos process.
The nullspace variant of the constraint-preconditioned version, Lanczos CP-CG, generates iterates \(\widehat{x}_k\)
so as to minimize the energy norm of the error, i.e.,
\[
  \|\widehat{e}_k\|_{\widehat{M}}^2  = \widehat{e}_k^T \widehat{M} \widehat{e}_k,
\]
where \(\widehat{e}_k = \widehat{x}_k - \widehat{x}_*\), and \(\widehat{x}_*\) is the exact solution of~\eqref{eq:reduced-sys}.
The definitions~\eqref{eq:reduced-sys-2} yield
\begin{align*}
  \|\widehat{e}_k\|_{\widehat{M}}^2 & =
  (\widehat{x}_k - \widehat{x}_*)^T Z^T M Z (\widehat{x}_k - \widehat{x}_*)
  \\ & =
  (x_k - x_*)^T A (x_k - x_*) + (w_k - w_*)^T F^{-1} (w_k - w^*)
  \\ & =
  (x_k - x_*)^T A (x_k - x_*) + (y_k - y_*)^T C (y_k - y_*),
\end{align*}
where we used again the relationship \(w_k = -F E^T y_k\) between iterates of P-CG and CP-CG.
%%%%% previous text
% $(x_k, y_k)$ such that the error vector
% $e_k = [x_k - x_* \,;\, y_k - y_*]$, where $[x_* \,;\, y_*]$ is the solution of~\eqref{eq:rsp}, satisfies
% \begin{equation}
% \label{eq:proj-error}
%   e_k^T M e_k \,
%   = \, \left \| \, \widehat{x}_k - \widehat{x}_* \, \right \|_{\widehat{M}}^2
%   = \, \min \left \| \, \widehat{x} - \widehat{x}_* \,\right \|_{\widehat{M}}^2
%      \mbox{ s.t. } \widehat{x} \in \widehat{x}_0 + \widehat{\mathcal{K}}_k .
% \end{equation}
For Lanczos CP-CG to be applicable, $\widehat{M}$ must be positive definite, which occurs when the sum
of the number of negative eigenvalues of \(K\) and \(C\) is \(m\) \cite[Theorem~$2.1$]{dollar-gould-schilders-wathen-2006}.

We can derive a ``traditional'' CP-CG implementation by applying the usual transformations to the Lanczos CP-CG.
The result coincides with the implementation of \cite{dollar-gould-schilders-wathen-2006}, although the latter authors assume that \(B\) has full row rank for specific purposes.
It is also equivalent to that of \cite{cafieri-dapuzzo-desimone-diserafino-2007a} for~\eqref{eq:rsp} with positive definite $C$.
% \smarttodo{Why do Dollar and al. require that $B$ have full rank???
% Daniela: in Section~3 of their  paper, they exploit the fact that $B$ has a submatrix $B_1 \in \R^{m \times m}$ of rank $m$}
The above suggests that CP-CG is appropriate to solve saddle-point systems
in constrained optimization where~\eqref{eq:rsp} is used to minimize a quadratic model of
a penalty function and sufficient decrease of this quadratic model is sought, such as in trust-region methods.

Our last example considers SYMMLQ \citep{paige-saunders-1975}, which does not require
$\widehat{M}$ to be positive definite but, like CG, requires~\eqref{eq:rsp} to be consistent.
Its constraint-preconditioned version, CP-SYMMLQ, computes $[x_k \,;\, y_k]$ so as to minimize the error in a norm defined by the preconditioner, i.e.,
\begin{align*}
  \widehat{e}_k^T \widehat{P}^{-1} \widehat{e}_k & =
  (\widehat{x}_k - \widehat{x}_*)^T (Z^T P Z)^{-1} (\widehat{x}_k - \widehat{x}_*) \\ & =
  (\widehat{x}_k - \widehat{x}_*)^T Z^T Z (Z^T P Z)^{-1} Z^T Z (\widehat{x}_k - \widehat{x}_*) \\ & =
  \begin{bmatrix}
    x_k - x_* \\ w_k - w_*
  \end{bmatrix}^T
  P_G
  \begin{bmatrix}
    x_k - x_* \\ w_k - w_*
  \end{bmatrix}.
  % (x_k - x_*)^T G (x_k - x^*) + (y_k - y_*)^T C (y_k - y_*),
\end{align*}
where we used similar identifications as above and assumed, without loss of generality, that \(Z\) has orthonormal columns.
In other words, if we define
\begin{equation}
  \label{eq:ex-ew}
  \begin{bmatrix}
       G &        & B^T
    \\   & F^{-1} & E^T
    \\ B & E      &
  \end{bmatrix}
  \begin{bmatrix}
       e_x
    \\ e_w
    \\ \bar{e}
  \end{bmatrix}
  =
  \begin{bmatrix}
       x_k - x_*
    \\ w_k - w_*
    \\ 0
  \end{bmatrix},
\end{equation}
then
\[
  \widehat{e}_k^T \widehat{P}^{-1} \widehat{e}_k =
  (x_k - x_*)^T e_x + (w_k - w_*)^T e_w =
  e_x^T G e_x + e_w F^{-1} e_w.
\]
By~\eqref{eq:w} and~\eqref{eq:ex-ew}, there exists a vector \(e_y\) such that \(e_w = -F E^T e_y\), and thus \(e_w^T F^{-1} e_w = e_y^T C e_y\).
The second block row of~\eqref{eq:ex-ew} premultiplied by \(E\) yields \(EF (w_k - w_*) - C \bar{e} = -C e_y\), so that~\eqref{eq:ex-ew} can be written as
\[
  \begin{bmatrix}
    G & \phantom{-}B^T \\
    B & -C
  \end{bmatrix}
  \begin{bmatrix}
    e_x \\ e_y
  \end{bmatrix}
  =
  \begin{bmatrix}
    x_k - x_*
    \\ 0
  \end{bmatrix}.
\]
Finally, CP-SYMMLQ minimizes
\[
  \widehat{e}_k^T \widehat{P}^{-1} \widehat{e}_k = e_x^T G e_x + e_y^T C e_y.
\]

% $$
%     e_k^T P^{-1} e_k \,
%     = \, \left \| \, \widehat{x}_k - \widehat{x}_* \, \right \|_{\widehat{P}}^2
%     = \, \min \left \| \, \widehat{x} - \widehat{x}_* \, \right \|_{\widehat{P}}^2
%     \mbox{ s.t. } \widehat{x} \in \widehat{x}_0 + \widehat{M} \, \widehat{\mathcal{K}}_{k-1}.
% $$
% \smarttodo{Check the minimization properties of the CP-Krylov solvers.}

\section{Constraint-Preconditioned Arnoldi Process and Associated Krylov
Solvers\label{sec:cp-arnoldi-and-methods}}

A constraint-preconditioned version of the Arnoldi process can be derived
by reasoning as in \Cref{sec:cp-lanczos}, obtaining \Cref{alg:cp-arnoldi}.
The equivalence between the projected version (\Cref{alg:proj-arnoldi} in~\Cref{sec:processes})
and the constraint-preconditioned version is stated in \Cref{thm:cp-arnoldi-equiv}, which is akin
to \Cref{thm:cp-lanczos-equiv}.

\begin{algorithm}[ht]
    \caption{Constraint-Preconditioned Arnoldi Process\label{alg:cp-arnoldi}}
                                                       % for \eqref{eq:rsp}}
    \begin{algorithmic}[1]
      \State choose $[x_0 \,;\, q_0]$ such that \(B x_0 - C q_0 = 0\) \Comment{initial guess}
      \State $p_0 = 0$ \Comment{initial Arnoldi vector}
      \State $u_0 = b - A x_0$, \ $t_0 = C q_0$
      \State \label{alg:cpa-line-proj2-0}%
        $[\bar{p}_1 \,;\, \bar{z}_1] \leftarrow$ solution of~\eqref{eq:proj-2} with right-hand side $[u_0\,;\, -t_0]$
      \State $p_1 = \bar{p}_1$                          \label{alg:cpa-line-p1}
      \State $q_1 = q_0 - \bar{z}_1$                    \label{alg:cpa-line-q1}
      \State $h_{1,0} = (p_1^T u_0 + q_1^T t_0)^{\tfrac{1}{2}}$
%      \IfOneLine {$h_{1,0} \ne 0$} \label{alg:cpa-line-scale1} \ $p_1 = p_1 / h_{1,0}$, \ $q_1 = q_1 / h_{1,0}$
      \If {$h_{1,0} \ne 0$}
           \State \label{alg:cpa-line-scale1} \ $p_1 = p_1 / h_{1,0}$, \ $q_1 = q_1 / h_{1,0}$
      \EndIf
      \State $k = 1$
      \While {$h_{k,k-1} \ne 0$}
         \State $u_k = A p_k$, \ $t_k = C q_k$
         \State \label{alg:cpa-line-proj2-k}
            $[\bar{p}_{k+1}\,;\, \bar{z}_{k+1}] \leftarrow$ solution of~\eqref{eq:proj-2} with right-hand side $[u_k\,;\, -t_k]$
         \State $p_{k+1} = \bar{p}_{k+1}$
         \State $q_{k+1} = q_k - \bar{z}_{k+1}$
         \For {$i = 1,\ldots,k$}
              \State $h_{i,k} = p_i^T u_k + q_i^T t_k$ \Comment{$= p_i^T A p_k + q_i^T C q_k$}
              \State $p_{k+1} = p_{k+1} - h_{i,k} p_i$
              \State $q_{k+1} = q_{k+1} - h_{i,k} q_i$
         \EndFor
         \State $h_{k+1,k} = (p_{k+1}^T u_k + q_{k+1}^T t_k)^{\tfrac{1}{2}}$ \Comment{$= (p_{k+1}^T A p_k + q_{k+1}^T C q_k)^{\tfrac{1}{2}}$}
%        \IfOneLine {$h_{k,k+1} \ne 0$} \label{alg:cpa-line-scalek+1} \ $p_{k+1} = p_{k+1} / h_{k,k+1}$, \ $q_{k+1} = q_{k+1} / h_{k,k+1}$
         \If {$h_{k+1,k} \ne 0$}
             \State \label{alg:cpa-line-scalek+1} \ $p_{k+1} = p_{k+1} / h_{k+1,k}$, \ $q_{k+1} = q_{k+1} / h_{k+1,k}$
        \EndIf
         \State $k = k+1$
       \EndWhile
    \end{algorithmic}
\end{algorithm}

% \par
% \begin{theorem} \label{th:equiv-cp-arnoldi}
% Let $E$ and $F$ be the matrices defined in~\eqref{eq:factc}.
% Then, \Cref{alg:proj-arnoldi} with starting guess $[x_0\,;\, w_0]$ is equivalent
% to \Cref{alg:cp-arnoldi} with starting guess $[x_0\,;\, q_0]$ such that
% $w_0 = -FE^T q_0$. In particular, for all~$k$, the vectors $v_{k,x}$ and $v_{k,w}$
% and the scalars $h_{i,k}$ $(i=1,\ldots,k+1)$ in \Cref{alg:proj-arnoldi} are equal
% to the vectors $p_k$ and $FE^T q_k$ and to the scalars $h_{i,k}$
% in \Cref{alg:cp-arnoldi}, respectively.
% \end{theorem}
\begin{shadytheorem}
  \label{thm:cp-arnoldi-equiv}
  Let $E$ and $F$ be as defined in~\eqref{eq:factc} and $G$ chosen to satisfy \Cref{assum:Phat-spd}.
  Let $q_0 \in \R^m$ be arbitrary.
  Then, \Cref{alg:proj-arnoldi} in~\Cref{sec:processes} with starting guesses $x_0 \in \R^m$ and $w_0 = -FE^T q_0$,
  such that $B x_0 + E w_0 = 0$, is equivalent to \Cref{alg:cp-arnoldi} with starting guesses $x_0$ and $q_0$.
  In particular, for all~$k$, the vectors $v_{k,x}$ and $v_{k,w}$ and the scalars $h_{i,k}$ in \Cref{alg:proj-arnoldi}
  are equal to the vectors $p_k$ and $FE^T q_k$, and to the scalars $h_{i,k}$ in \Cref{alg:cp-arnoldi}, respectively.
\end{shadytheorem}

As in the case of the Lanczos process, the CP-Arnoldi process is equivalent to applying
the corresponding standard Arnoldi process to system~\eqref{eq:rsp} with preconditioner~\eqref{eq:cp}
(see \Cref{alg:cp-arnoldi-full-space} in \Cref{sec:processes}), as stated in the next theorem.
\begin{shadytheorem}
  \label{thm:cp-arnoldi-equiv-2}
  Let \Cref{alg:cp-arnoldi} be initialized with \(x_0 \in \R^n\) such that \(B x_0 = 0\), \(q_0 = 0 \in \R^m\), and \Cref{alg:cp-arnoldi-full-space}
  be initialized with the same \(x_0\) and \(y_0 \in \R^m\) such that~\eqref{eq:2ndeq0_x0y0} is satisfied.
  Then, for all \(k \geq 0\), \(v_{k,x} = p_k\) and \(v_{k,y} = -q_k\), where \([v_{k,x}\,;\, v_{k,y}]\) is the \(k\)-th Arnoldi vector
  generated in \Cref{alg:cp-arnoldi-full-space}, and \(p_k\) and \(q_k\) are the \(k\)-th Arnoldi vectors generated in \Cref{alg:cp-arnoldi}.
  In addition, the scalars \(h_{i,k}\) computed at each iteration are the same in both algorithms.
\end{shadytheorem}

\Cref{thm:cp-arnoldi-equiv-2} allows us to develop a constraint-preconditioned variant of any Krylov method
based on the Arnoldi process, using a starting guess satisfying~\eqref{eq:2ndeq0_x0y0}.
Such variants are equivalent to their standard counterparts preconditioned with~\eqref{eq:cp},
but are computationally cheaper, as in the case of Lanczos-based methods.
Furthermore, CP-Krylov versions of optimal Arnoldi-based Krylov methods preserve the minimization
properties of these methods in the sense explained in \Cref{sec:cp-l-methods}.
For example, the constraint-preconditioned version of GMRES \citep{saad-schultz-1986} minimizes
the norm of the deviation of the residual from \(\Range(N)\) similarly to MINRES.
Below, we denote GMRES(\(\ell\)) the variant of GMRES that is restarted every \(\ell\) iterations.

Obtaining constraint-preconditioned versions of GMRES($\ell$) and DQGMRES is straighforward,
by restarting and truncating the CP-Arnoldi basis generation process, respectively, as in the standard case~\citep{saad-2003}.
Note that DQGMRES with memory~$2$, i.e., with orthogonalization of each Arnoldi vector against the two previous
vectors only, is equivalent to CP-MINRES in exact arithmetic when $A$ is symmetric.
In finite precision arithmetic, DQGMRES with a larger memory may dampen the loss of orthogonality among
the Lanczos vectors and act as a local reorthogonalization procedure, although we did not observe significant differences in \Cref{sec:experiments}.

\citet[Theorems~\(4.1\) and~\(4.3\)]{dollar-2007} establishes that if \(C\) is positive semi-definite of rank \(p\),
\(P^{-1} K\) has an eigenvalue at~\(1\) of multiplicity \(2m - p\), while the remaining \(n - m + p\) eigenvalues are
defined by a generalized eigenvalue problem. A remark after \cite[Theorem~\(4.1\)]{dollar-2007} states that
\Cref{assum:Phat-spd} ensures that all eigenvalues are real. In addition, the dimension of the Krylov space is at most \(\min(n - m + p + 2, \, n + m)\).
Inspection reveals that \citeauthor{dollar-2007}'s proofs of those results do not use the fact that \(A\) is symmetric; the results hold for general \(A\).
\cite{loghin-2017}
% \smarttodo{Modified a bit the note about Loghin. OK?}
establishes similar results on the eigenvalues of non-regularized saddle-point matrices for general~\(A\) and general~\(G\).
Clustering eigenvalues accelerates convergence of nonsymmetric Krylov solvers in many practical cases, although the convergence behavior
of such solvers is not fully characterized by the eigenvalues \citep{greenbaum-ptak-strakos-1996}.
% \smarttodo{I was probably unclear in my question about \citep{loghin-2017}. I wondered if we had to cite this paper because there it is
% observed that the well-known results about the eigenvalues of constraint-preconditioned non-regularized saddle-point systems also hold
% for general \(A\). In any case, I don't think \citep{loghin-2017} must be cited together with \citep{greenbaum-ptak-strakos-1996}. Maybe
% we can remove this citation. DO YOU AGREE?}

%\section{The Constraint-Preconditioned Orthogonal Tridiagonalization Process\label{sec:ssy}}
%\citep{saunders-simon-yip-1988}
%
%USYMQR, USYMLQ

\section{Implementation Issues}\label{sec:implementation}

%Our MATLAB implementation provides the \texttt{cpkrylov} library, containing the constraint-preconditioned
%variants of the Lanczos-CG, MINRES, SYMMLQ, GMRES($\ell$) and DQGMRES methods for~\eqref{eq:rsp}.
We implemented the constraint-preconditioned variants of the Lanczos-CG, MINRES, SYMMLQ, GMRES($\ell$)
and DQGMRES methods for~\eqref{eq:rsp} in a MATLAB library named \texttt{cpkrylov}.
For completeness, we also included in the library an implementation of the CP-CG method in the form
given by \cite{dollar-gould-schilders-wathen-2006}. We think that \texttt{cpkrylov} can be useful as a basis
for the development of more sophisticated numerical software.

All solvers are accessed via a common interface exposed by the main driver \texttt{reg\_cpkrylov()}, which
performs pre-processing operations, calls the requested solver, performs post-processing operations, and returns solutions and statistics to the user.
\texttt{cpkrylov} is freely available from \https{github.com/optimizers/cpkrylov}.

Because \(A\) is never required as an explicit matrix, we allow the user to supply it as an abstract linear operator as implemented in the Spot linear operator toolbox\footnote{\http{www.cs.ubc.ca/labs/scl/spot}}.
Spot allows us to use the familiar matrix notation with operators for which a representation as an explicit matrix is unavailable or inefficient.
This affords the user flexibility in defining \(A\) while keeping the implementation of the various Krylov methods as readable as if \(A\) were a matrix.
%\smarttodo[inline]{Daniela: any reference to spot must be removed if the spot-based implementation is inefficient (see also the next paragraph).}
%\smarttodo[inline,color=white,linecolor=gray]{Dom: did the updated files I sent you not correct the issue?}
%\smarttodo[inline]{Problem solved!}

% The constraint-preconditioned versions of the Lanczos CG, MINRES,
% SYMMLQ, GMRES($\ell$) and DQGMRES methods for solving system~\eqref{eq:rsp}
% were implemented in MATLAB, not only to allow an easy use of
% these solvers, but also to make available ``templates'' that can be a basis for
% the development of more sophisticated software. In order to obtain more flexible implementations
% and simplify code development, the Spot linear operator toolbox was exploited
% (\url{http://www.cs.ubc.ca/labs/scl/spot/}), which extends the MATLAB built-in matrix notation to linear
% operators for which a representation as an explicit matrix is unavailable.
% All the solvers can be accessed through the same interface, i.e., by calling a driver function
% named \verb|reg_cpkrylov|, which in turn calls the CP-Krylov solver specified by the user and performs some
% pre-processing and post-processing operations, as detailed next. For completeness,
% an implementation of the CP-CG solver by~\cite{dollar-gould-schilders-wathen-2006} was
% developed too. The MATLAB codes are freely available from
% \url{https://github.com/optimizers/cpkrylov}.

\cite{gould-hribar-nocedal-2001,gould-orban-rees-2014} observe that the numerical stability of projected Krylov solvers depends on keeping $[x_k\,;\, w_k]$ in \(\Null(N)\).
While the iterates lie in the nullspace in exact arithmetic, $[x_k\,;\, w_k]$ may have a non-negligible component in $\Range(N^T)$ because of roundoff error.
In turn, the stability of CP-Krylov solvers depends on how accurately \([x_k\,;\, y_k]\) satisfies
$$
   B x_k - C y_k = 0.
$$
\cite{gould-hribar-nocedal-2001} suggest to increase the accuracy by applying iterative refinement after solving~\eqref{eq:proj-2} with a direct method.
In \texttt{cpkrylov}, the constraint preconditioner $P$
is implemented as a Spot operator \texttt{P}
such that writing \texttt{z = P*r}, where \texttt{z = [z1 ; z2]} and \texttt{r = [r1 ; r2]}, corresponds to solving
\begin{equation}
  \label{eq:applP}
   \begin{bmatrix}
   G  &  \phantom{-}B^T \\
   B  &  -C\phantom{^T}
   \end{bmatrix}
   \begin{bmatrix}
   z_1 \\ z_2
   \end{bmatrix}
   =
   \begin{bmatrix}
   r_1 \\ r_2
   \end{bmatrix},
\end{equation}
and automatically performing iterative refinement if requested by the user or if the residual norm of~\eqref{eq:applP} exceeds a given tolerance.

An alternative approach to minimizing the size of the component of \([x_k\, ; w_k]\) in \(\Range(N^T)\) suggested by \cite{gould-hribar-nocedal-2001} is to perform \emph{iterative semi-refinement}.
The latter consists in noting that the solution of~\eqref{eq:proj} is not affected (in exact arithmetic) if we add a vector lying in \(\Range(N^T)\) to \([u_x \, ; \, u_w]\) in the right-hand side.
Such a vector is available cheaply in the form of \([B^T \bar{z} \, ; \, E^T \bar{z}]\) where \(\bar{z}\) is the trailing segment of the solution of the most recent projection step~\eqref{eq:proj}, and \(\bar{z} = 0\) at the first projection step.
The net result is that instead of~\eqref{eq:proj-2}, we solve
\begin{equation}
  \label{eq:proj-2-isr}
  \begin{bmatrix}
   G & \phantom{-}B^T \\
   B & -C\phantom{^T}
  \end{bmatrix}
  \begin{bmatrix}
    \bar{p}_{k+1} \\ \bar{z}_{k+1}
  \end{bmatrix}
  =
  \begin{bmatrix}
    u_{k,x} - B^T \bar{z}_k \\ -(t_k - C \bar{z}_k)
  \end{bmatrix},
  \quad \bar{z}_0 := 0.
\end{equation}

By default, the matrix of~\eqref{eq:applP} is factorized by way of MATLAB's \texttt{ldl()}.
Spot allows us to separate the implementation of the preconditioner from that of other phases of solvers, so that future extensions
to the former (e.g., the case where applying \(P\)
results from a different factorization) will not require changes to the latter.
Our implementation of $P$ is transparent to the user, who must only pass the matrices $G$, $B$ and $C$ to \texttt{reg\_cpkrylov()}.

All CP-Krylov solvers stop when
\begin{equation} \label{eq:stop-crit}
  \|r_{\text{P},k}\|_{[P]} \le \epsilon_a + \|r_{\text{P},0}\|_{[P]} \, \epsilon_r,
\end{equation}
%\smarttodo{Daniela: check the stopping criterion.}
where \(\|r_{\text{P},k}\|_{[P]}\) is defined in~\eqref{eq:norm-dev} (or, equivalently, in~\eqref{eq:norm-dev-2}),
and $\epsilon_a$ and $\epsilon_r$ are tolerances given by the user (default values are also set in our
implementations). Note that, for all the CP-Krylov solvers except CP-DQGMRES, \(\|r_{\text{P},k}\|_{[P]}\) is obtained
as a byproduct of other computations performed in algorithm. CP-DQGMRES computes an estimate
of the residual norm only.
%, but this is considered reasonably good because the Arnoldi vectors are nearly orthogonal in an appropriate space.
A computationally cheap overestimate of the residual norm could be
used in the stopping criterion, but this may unnecessarily increase the number of iterations~\citep[Section~3.1]{saad-wu-1996}.
A maximum number of iterations can be also specified for all solvers.

%CP-SYMMLQ computes both the so-called CG and the LQ residuals (see Paige \& Saunders 1975),
%but the CG residual is used in the stopping criterion. The best solution among the CG and the LQ one
%is provided in output. The MINRES residual is also computed for comparison.

So far, we have considered the case where the last $m$ entries of the right-hand side of~\eqref{eq:rsp} are zero.
When the right-hand side has the general form $[b_1\,;\, b_2]$ with $b_2 \ne 0$, we can compute $\Delta x$ and $\Delta y$ such that
\begin{equation}
  \label{eq:BDx-CDy}
  B \Delta x - C \Delta y = b_2,
\end{equation}
by applying \(P\) to \([0 \,;\, b_2]\), and subsequently solve~\eqref{eq:rsp} with $b = b_1 - A \Delta x - B^T \Delta y$.
The solution of the original system is $[x + \Delta x\,;\, y + \Delta y]$.
These pre- and post-processing steps are implemented in \texttt{reg\_cpkrylov()}.

%By numerical experiments, we observed that the approximation $\bar{y}$ of $y_*$ computed
%directly by a CP-Krylov solver is generally comparable or less accurate than the approximation obtained
%with the following procedure. System~\eqref{eq:rsp} is equivalent to
%$$
%  \begin{bmatrix}
%    A+G & \phantom{-}B^T \\
%    B & -C\phantom{^T}
%  \end{bmatrix}
%  \begin{bmatrix}
%  x \\
%  y
%  \end{bmatrix}
%  =
%  \begin{bmatrix}
%  b + G x \\
%  0
%  \end{bmatrix}
%$$
%and hence $y_*$ could be obtained by solving
%\begin{equation}
%\label{eq:comp-y}
%  \begin{bmatrix}
%    G & \phantom{-}B^T \\
%    B & -C\phantom{^T}
%  \end{bmatrix}
%  \begin{bmatrix}
%  x \\
%  y
%  \end{bmatrix}
%  =
%  \begin{bmatrix}
%  b - A x_* + G x_*\\
%  0
%  \end{bmatrix},
%\end{equation}
%which is obviously impractical. Since $x_*$ is not available, we solve system~\eqref{eq:comp-y}
%with $\bar{x}$ instead of $x_*$ in the right-hand side, and take as approximation of $y_*$
%the vector $\underline{y}$ consisting of last $m$ entries of the solution.
%This corresponds to one more application of the constraint preconditioner and is performed
%by \verb|reg_cpkrylov| after running the Krylov solver. \smarttodo{Is this really true or does it depend
%on the way the random test problems are generated? How do we justify it? Can we simply say
%that this deserves further investigation, which is beyond the scope of this work?}

\section{Numerical Experiments\label{sec:experiments}}

We report results obtained by applying some solvers from the \texttt{cpkrylov} library to regularized saddle-point systems
arising in the application of the primal-dual interior point solver PDCO to convex quadratic programming problems
(see \https{web.stanford.edu/group/SOL/software/pdco/}).
PDCO solves linearly constrained optimization problems with a smooth convex objective function in the form
% $$
% \begin{array}{cl}
% \displaystyle \minimize{x \in \R^n}  & f(x) \\
% \mbox{subject to}                                                & Bx = c \\[3pt]
%                                                                           & l \le x \le u,
% \end{array}
% $$
% where $f: \R^n \to \R$ and $B \in \R^{m \times n}$.
% Actually, to ensure unique primal and dual solutions and improve stability, PDCO solves the problem
\begin{equation} \label{eq:reg_optim_pb}
\begin{array}{cl}
\displaystyle \minimize{x \in \R^n, \, r \in \R^m}  & f(x) + \frac12 \|D_1 x\|^2 + \frac12 \|r\|^2 \\
\mbox{subject to}                                                & Bx + D_2 r = c \\[3pt]
                                                                           & l \le x \le u,
\end{array}
\end{equation}
where $f: \R^n \to \R$ is smooth and convex, $B \in \R^{m \times n}$, and \(D_1\) and \(D_2\) are positive-definite diagonal matrices
that provide primal and dual regularization. In particular,
\(D_2\) determines whether \(Bx = c\) should be satisfied accurately or in the least-squares sense.

At each iteration of PDCO, a Newton step is applied to suitably perturbed KKT conditions associated with~\eqref{eq:reg_optim_pb}.
The Newton step requires the solution of a linear system, which can be cast into the form~\eqref{eq:rsp} by a combination of permutation
operations and/or inexpensive block eliminations. Possibly the most common saddle-point formulation is
$$
K_2 =
\begin{bmatrix}
   A & \phantom{-}B^T \\
   B & -C\phantom{^T}
\end{bmatrix}
=
\begin{bmatrix}
  H + D_1^2 + X_1^{-1}Z_1 + X_2^{-1}Z_2  & \phantom{-}B^T \\
  B                                                                 & -D_2^2
\end{bmatrix}
$$
where $H$ is the Hessian of the objective function at the current approximation of the optimal solution, $X_1 = \diag(x_1)$,  $X_2 = \diag(x_2)$,
$Z_1 = \diag(z_1)$, $Z_2 = \diag(z_2)$, $x_1 = x - l > 0$, $x_2 = u - x > 0$, and $z_1 > 0$ and $z_2 > 0$ are the corresponding dual variable estimates.
In our experiments, $H$ is constant because \(f\) is quadratic.
More details are available from \https{web.stanford.edu/group/SOL/software/pdco/pdco.pdf}.

Recently, unreduced KKT systems have attracted the interest of researchers because of their better spectral properties, especially
as the interior point iterates approach the solution of the optimization problem~\citep{greif-moulding-orban-2014,morini-simoncini-tani-2016}.
Other symmetric and unsymmetric saddle-point formulations are obtained with simple operations.
In particular, within PDCO we used the unreduced symmetric saddle-point formulation
%\smarttodo{I suggest \(K_{3.5}\) for consistency with other papers. Should \(X\) be \(-X\)?}
$$
K_{3.5} =
\begin{bmatrix}
    A & \phantom{-}B^T \\
    B & -C\phantom{^T}
\end{bmatrix}
=
\left[ \begin{array}{c|cc}
   H + D_1^2                  & \phantom{-}B^T                   & \phantom{-}Z^\frac12                    \\ \hline
   B\phantom{^\frac12}   & -D_2^2                                & \phantom{-}0\phantom{^\frac12}  \\
   Z^\frac12                     & \phantom{-}0\phantom{^T} & -X\phantom{^\frac12}                    \\
\end{array}  \right] ,
$$
where $X = \diag([x_1 \,; \, x_2])$ and $Z = \diag([z_1 \,; \, z_2])$.
We also considered the unsymmetric saddle-point formulation
%\smarttodo{I suggest \(K_{3p}\)}
$$
K_{3p} =
\begin{bmatrix}
    A & \phantom{-}B^T \\
    B & -C\phantom{^T}
\end{bmatrix}
=
\left[ \begin{array}{cc|c}
   H + D_1^2        & I   & \phantom{-}B^T                   \\
   -Z                     & X  & \phantom{-}0\phantom{^T}  \\ \hline
   \phantom{-}B   &  0   & -D_2^2                               \\
\end{array}  \right] ,
$$
which has the same structure as the saddle-point matrix in equation~(2.5a) of \citep{greif-moulding-orban-2014} up to a permutation.

For all the saddle-point formulations, the constraint preconditioner $P$ in~\eqref{eq:cp} was defined by choosing $G$ equal to the diagonal
of the leading block $A$. This is a common choice in interior point methods---see, e.g., \citep{dapuzzo-desimone-diserafino-2010}.
In our experiments, iterative refinement never needed to be performed.

The CP-Krylov solvers were stopped using an adaptive criterion that relates the accuracy in the solution of the KKT linear system
to the duality measure at the current interior point iteration, as suggested by \cite{cafieri-dapuzzo-desimone-diserafino-2007b}.
Thus, criterion~\eqref{eq:stop-crit} was applied by setting $\epsilon_r = 0$ and
$$
\epsilon_a = \max \left\{ \min \left\{ 10^{-2} \mu, 10^{-2} \right\}, \, 10^{-6} \right\},
$$
where $\mu$ is the the barrier parameter in PDCO. % duality measure.
% \smarttodo{rather the barrier parameter in PDCO, no?}

We run PDCO on the \textsf{CUTEst} \citep{gould-orban-toint-2015} problems reported in \Cref{tab:qp-prob}.
We use the models translated\footnote{\https{github.com/mpf/Optimization-Test-Problems}} into the AMPL modeling language \citep{AmplBook2}.
Our version of PDCO has been modified to take an optimization problem in the form of an instance of the \texttt{nlpmodel} class as argument, which is defined in the \texttt{model} Matlab package.\footnote{\https{github.com/optimizers/model}}
The \texttt{amplmodel} subclass of \texttt{nlpmodel} reads an AMPL \texttt{nl} file by way of the \texttt{AmplMEXInterface} package\footnote{\https{github.com/optimizers/AmplMexInterface}} and conforms to the \texttt{model} interface expected by our version of PDCO.
The rest of PDCO is identical to the original version.
For the problem received by PDCO to have the form~\eqref{eq:reg_optim_pb}, it is necessary to introduce slack variables.
In our implementation, linear inequalities \(\ell \leq Ax \leq u\) are transformed to \(Ax - s = 0\) and \(\ell \leq s \leq u\).
The transformation is performed by the \texttt{slackmodel} class, which receives an arbitrary instance of \texttt{nlpmodel}, including arbitrary instances of \texttt{amplmodel}, and adds slack variables as just described.
The options passed to PDCO, including scaling parameters, are the same as those described by \cite{orban-2015}.

The problems are chosen so that $H$ is non-diagonal; otherwise, for $K_2$ and $K_{3.5}$ the constraint
preconditioner would be equal to the saddle-point matrix. The table also shows the
% number of variables and linear constraint, as well as the
sizes of $K_2$ and $K_{3.5}$ (or $K_{3p}$).

\begin{table}[t!]
{\scriptsize
\begin{center}
\begin{tabular}{lrr}
\hline
Problem     & $K_2$ size                 & $K_{3.5} \; (K_{3p})$ size \\ \hline
cvxqp1\_s  &          $ 200 $  &       $ 400 $  \\
cvxqp1\_m & $     2000 $  &    $ 4000 $  \\
cvxqp1\_l   &  $ 20000 $ & $ 40000 $  \\
cvxqp2\_s  &         $ 150 $  &         $ 350 $  \\
cvxqp2\_m &     $ 1500 $  &     $ 3500 $  \\
cvxqp2\_l   &  $ 15000 $ & $ 35000 $  \\
cvxqp3\_s  &          $ 250 $  &         $ 450 $  \\
cvxqp3\_m &      $ 2500 $  &     $ 4500 $  \\
cvxqp3\_l   & $ 25000 $  &  $ 45000 $  \\
gouldqp3   &     $ 1397 $  &      $ 2795 $  \\
gouldqp2   &     $ 1397 $  &      $ 2795 $  \\
mosarqp1  &     $ 3900 $  &      $ 7100 $  \\
mosarqp2  &     $ 3900 $  &      $ 3600 $  \\
stcqp1       &     $ 8201 $  &    $ 16395 $  \\
stcqp2       &     $ 8201 $  &    $ 16395 $  \\ \hline
\end{tabular}
\end{center}
}
\par \vspace{3pt}
\caption{\textsf{CUTEst} problems used in the experiments.\label{tab:qp-prob}}
\end{table}

We ran PDCO with the saddle-point matrices $K_2$ and $K_{3.5}$, using CP-CG, CP-MINRES,
CP-DQGMRES($\ell$) and CP-GMRES($\ell$) as Krylov solvers. By CP-DQGMRES($\ell$) we denote CP-DQGMRES with
memory parameter $\ell$, i.e., the number of Arnoldi vectors to be stored in the truncated CP-Arnoldi process.
We set $\ell = 2$; in this case CP-DQGMRES is equivalent to CP-MINRES in exact arithmetic.
We also ran PDCO with $K_{3p}$ using CP-DQGMRES($\ell$) and CP-GMRES($\ell$) with various values of $\ell$.
The goal of the experiments is to illustrate the behavior of CP-Krylov solvers inside an interior-point method.

PDCO was run on a 2.5 GHz Intel Core i7 processor with 16 GB of RAM,
4~MB of L3 cache and the macOS 10.13.6 operating system, using MATLAB R2018b.
Execution times were measured in seconds, by using the MATLAB function \texttt{timeit},
which removes some of the noise inherent to time measurements by calling a specified function
multiple times and returning the median of the measurements.

\begin{table}[h!]
{\scriptsize
\begin{center}
\begin{tabular}{lrrrrr}
\hline
name & outer it & inner it & PDCO time & prec time & solve time \\
\hline
      cvxqp1\_s &   17 &     80 & \( 1.025\)e\(-01\) & \( 4.137\)e\(-02\) & \( 4.276\)e\(-02\) \\
      cvxqp1\_m &   19 &    103 & \( 2.233\)e\(-01\) & \( 7.196\)e\(-02\) & \( 1.190\)e\(-01\) \\
      cvxqp1\_l &   20 &    138 & \( 1.367\)e\(+00\) & \( 4.195\)e\(-01\) & \( 6.623\)e\(-01\) \\
      cvxqp2\_s &   17 &     80 & \( 6.875\)e\(-02\) & \( 3.357\)e\(-02\) & \( 2.102\)e\(-02\) \\
      cvxqp2\_m &   19 &    118 & \( 1.775\)e\(-01\) & \( 4.506\)e\(-02\) & \( 1.032\)e\(-01\) \\
      cvxqp2\_l &   20 &    140 & \( 9.295\)e\(-01\) & \( 1.468\)e\(-01\) & \( 5.485\)e\(-01\) \\
      cvxqp3\_s &   20 &     72 & \( 8.294\)e\(-02\) & \( 4.222\)e\(-02\) & \( 2.336\)e\(-02\) \\
      cvxqp3\_m &   19 &     99 & \( 2.731\)e\(-01\) & \( 1.098\)e\(-01\) & \( 1.309\)e\(-01\) \\
      cvxqp3\_l &   20 &    137 & \( 1.236\)e\(+00\) & \( 5.764\)e\(-01\) & \( 3.727\)e\(-01\) \\
       gouldqp3 &   10 &     20 & \( 5.697\)e\(-02\) & \( 2.671\)e\(-02\) & \( 1.619\)e\(-02\) \\
       gouldqp2 &   11 &     26 & \( 7.421\)e\(-02\) & \( 3.219\)e\(-02\) & \( 2.546\)e\(-02\) \\
       mosarqp1 &   17 &     54 & \( 2.962\)e\(-01\) & \( 9.205\)e\(-02\) & \( 1.595\)e\(-01\) \\
       mosarqp2 &   16 &     73 & \( 2.344\)e\(-01\) & \( 8.506\)e\(-02\) & \( 1.212\)e\(-01\) \\
         stcqp1 &   15 &    124 & \( 4.959\)e\(+00\) & \( 3.559\)e\(+00\) & \( 1.219\)e\(+00\) \\
         stcqp2 &   16 &    199 & \( 7.196\)e\(-01\) & \( 1.006\)e\(-01\) & \( 5.132\)e\(-01\) \\
\hline
\end{tabular}
\end{center}
}
\par \vspace{3pt}
\caption{Results for $K_2$ with solver CP-CG. Times are in seconds.\label{tab:res-k2-cpcg}}
\end{table}
\begin{table}[h!]
{\scriptsize
\begin{center}
\begin{tabular}{lrrrrr}
\hline
name & outer it & inner it & PDCO time & prec time & solve time \\
\hline
      cvxqp1\_s &   17 &     80 & \( 1.134\)e\(-01\) & \( 4.936\)e\(-02\) & \( 3.212\)e\(-02\) \\
      cvxqp1\_m &   19 &    103 & \( 2.231\)e\(-01\) & \( 7.145\)e\(-02\) & \( 1.182\)e\(-01\) \\
      cvxqp1\_l &   20 &    137 & \( 1.312\)e\(+00\) & \( 4.185\)e\(-01\) & \( 6.098\)e\(-01\) \\
      cvxqp2\_s &   17 &     80 & \( 6.928\)e\(-02\) & \( 3.350\)e\(-02\) & \( 2.168\)e\(-02\) \\
      cvxqp2\_m &   19 &    118 & \( 1.870\)e\(-01\) & \( 4.766\)e\(-02\) & \( 1.073\)e\(-01\) \\
      cvxqp2\_l &   20 &    140 & \( 9.056\)e\(-01\) & \( 1.506\)e\(-01\) & \( 5.127\)e\(-01\) \\
      cvxqp3\_s &   20 &     72 & \( 8.724\)e\(-02\) & \( 4.401\)e\(-02\) & \( 2.511\)e\(-02\) \\
      cvxqp3\_m &   19 &     99 & \( 2.869\)e\(-01\) & \( 1.164\)e\(-01\) & \( 1.356\)e\(-01\) \\
      cvxqp3\_l &   20 &    136 & \( 1.236\)e\(+00\) & \( 5.739\)e\(-01\) & \( 3.767\)e\(-01\) \\
       gouldqp3 &   10 &     20 & \( 5.260\)e\(-02\) & \( 2.453\)e\(-02\) & \( 1.531\)e\(-02\) \\
       gouldqp2 &   11 &     26 & \( 6.434\)e\(-02\) & \( 2.757\)e\(-02\) & \( 2.210\)e\(-02\) \\
       mosarqp1 &   17 &     54 & \( 3.383\)e\(-01\) & \( 9.735\)e\(-02\) & \( 1.671\)e\(-01\) \\
       mosarqp2 &   16 &     73 & \( 2.340\)e\(-01\) & \( 8.128\)e\(-02\) & \( 1.246\)e\(-01\) \\
         stcqp1 &   15 &    124 & \( 5.115\)e\(+00\) & \( 3.624\)e\(+00\) & \( 1.307\)e\(+00\) \\
         stcqp2 &   16 &    196 & \( 7.310\)e\(-01\) & \( 1.017\)e\(-01\) & \( 5.215\)e\(-01\) \\
\hline
\end{tabular}
\end{center}
}
\par \vspace{3pt}
\caption{Results for $K_2$ with solver CP-MINRES. Times are in seconds.\label{tab:res-k2-minres}}
\end{table}

\Cref{tab:res-k2-cpcg,tab:res-k2-minres,tab:res-k3.5-cg,tab:res-k3.5-minres} summarize the results obtained with $K_2$ and $K_{3.5}$
using CP-CG and CP-MINRES.
For each problem, \replaced[comment=R1.1]{``outer it'' is the number of outer interior-point iterations, ``inner it'' is the cumulative number of inner Krylov iterations, ``PDCO time'' is the total run time reported by PDCO, and ``prec time'' and ``solve time'' are the cumulative times to assemble and factorize the constraint preconditioner and to solve the linear systems, respectively.}{the number of PDCO and CP-Krylov iterations is reported, as well as the time for running PDCO, and the total time for building the constraint preconditioner and for solving the linear systems.}
%\smarttodo{Should we write something like ``Of course, small differences in the times for building the preconditioners for the same
%same saddle-point matrices within PDCO are due to the noise in the time measurements, which cannot be removed completely''? It's quite obvious ...}
We see that CP-MINRES performs a slightly smaller number of iterations than CP-CG on some problems,
which may be beneficial if very large systems are solved.
CP-MINRES is adequate in the context of a linesearch inexact Newton method such as PDCO because it reduces the residual norm monotonically by design.
\cite{fong-saunders-2012} observe that MINRES possesses other desirable properties that are generally attributed to CG.

We also observe that the number of CP-Krylov iterations with $K_{3.5}$ is always smaller than with $K_2$,
which may be due to the smaller condition number of $K_{3.5}$ \citep{greif-moulding-orban-2014,morini-simoncini-tani-2016}.
On cvxqp3\_s, $K_{3.5}$ also results in a smaller number of PDCO iterations.
\added[comment=R1.1]{On gouldqp2, \(K_{3.5}\) results in fewer inner iterations than outer iterations because the initial guess satisfies the stopping condition of the first five subproblems, resulting in zero inner iterations for those outer iterations. This behavior does not occur with \(K_2\), which produces different multiplier estimates.}
We used the MATLAB function \texttt{condest} to estimate the condition numbers of \(K_2\) and \(K_{3.5}\)
encountered during the PDCO iterations for each problem.
On the cvxqp problems, the largest value of \texttt{condest($K_{3.5}$)} is between three and four orders of magnitude smaller than the largest
value of \texttt{condest($K_2$)}.
The factor is between five and seven orders of magnitude on the gouldqp problems, one to two orders on the mosarqp problems,
and two to three orders on the stcqp problems. While such measurements do not tell the whole story and it would be more accurate to measure
the condition number of~\eqref{eq:reduced-sys-2}, they tend to confirm that the condition number of \(K_{3.5}\) is provably bounded if strict
complementarity is satisfied.
\replaced[comment=R2.4]{On our test set, the PDCO time reported for formulation \(K_2\) is almost always slightly smaller than that for \(K_{3.5}\).}{However, by looking at the times for PDCO we see that there is no clear winner between $K_2$ and $K_{3.5}$.}

\begin{table}[t!]
{\scriptsize
\begin{center}
\begin{tabular}{lrrrrr}
\hline
name & outer it & inner it & PDCO time & prec time & solve time \\
\hline
      cvxqp1\_s &   17 &     66 & \( 1.162\)e\(-01\) & \( 4.206\)e\(-02\) & \( 3.165\)e\(-02\) \\
      cvxqp1\_m &   19 &     86 & \( 3.685\)e\(-01\) & \( 9.864\)e\(-02\) & \( 2.275\)e\(-01\) \\
      cvxqp1\_l &   20 &    120 & \( 3.761\)e\(+00\) & \( 1.159\)e\(+00\) & \( 2.220\)e\(+00\) \\
      cvxqp2\_s &   17 &     64 & \( 7.814\)e\(-02\) & \( 4.010\)e\(-02\) & \( 2.098\)e\(-02\) \\
      cvxqp2\_m &   19 &    101 & \( 2.769\)e\(-01\) & \( 6.854\)e\(-02\) & \( 1.712\)e\(-01\) \\
      cvxqp2\_l &   20 &    123 & \( 2.483\)e\(+00\) & \( 4.038\)e\(-01\) & \( 1.753\)e\(+00\) \\
      cvxqp3\_s &   18 &     50 & \( 8.543\)e\(-02\) & \( 4.123\)e\(-02\) & \( 2.774\)e\(-02\) \\
      cvxqp3\_m &   19 &     81 & \( 4.584\)e\(-01\) & \( 1.793\)e\(-01\) & \( 2.384\)e\(-01\) \\
      cvxqp3\_l &   20 &    119 & \( 4.253\)e\(+00\) & \( 1.726\)e\(+00\) & \( 2.163\)e\(+00\) \\
       gouldqp3 &   10 &     12 & \( 5.470\)e\(-01\) & \( 4.725\)e\(-01\) & \( 5.119\)e\(-02\) \\
       gouldqp2 &    9 &      5 & \( 8.842\)e\(-02\) & \( 4.019\)e\(-02\) & \( 3.323\)e\(-02\) \\
       mosarqp1 &   17 &     39 & \( 5.908\)e\(-01\) & \( 2.126\)e\(-01\) & \( 3.048\)e\(-01\) \\
       mosarqp2 &   16 &     60 & \( 3.951\)e\(-01\) & \( 1.590\)e\(-01\) & \( 2.025\)e\(-01\) \\
         stcqp1 &   15 &    110 & \( 5.622\)e\(+00\) & \( 3.541\)e\(+00\) & \( 1.883\)e\(+00\) \\
         stcqp2 &   16 &    183 & \( 1.841\)e\(+00\) & \( 1.721\)e\(-01\) & \( 1.527\)e\(+00\) \\
\hline
\end{tabular}
\end{center}
}
\par \vspace{3pt}
\caption{Results for $K_{3.5}$ with solver CP-CG. Times are in seconds.\label{tab:res-k3.5-cg}}
\end{table}
\begin{table}[h!]
{\scriptsize
\begin{center}
\begin{tabular}{lrrrrr}
\hline
name & outer it & inner it & PDCO time & prec time & solve time \\
\hline
      cvxqp1\_s &   17 &     65 & \( 9.533\)e\(-02\) & \( 4.257\)e\(-02\) & \( 3.340\)e\(-02\) \\
      cvxqp1\_m &   19 &     86 & \( 3.882\)e\(-01\) & \( 1.030\)e\(-01\) & \( 2.462\)e\(-01\) \\
      cvxqp1\_l &   20 &    119 & \( 3.362\)e\(+00\) & \( 9.801\)e\(-01\) & \( 2.039\)e\(+00\) \\
      cvxqp2\_s &   17 &     64 & \( 6.759\)e\(-02\) & \( 3.480\)e\(-02\) & \( 1.855\)e\(-02\) \\
      cvxqp2\_m &   19 &    101 & \( 2.633\)e\(-01\) & \( 6.210\)e\(-02\) & \( 1.668\)e\(-01\) \\
      cvxqp2\_l &   20 &    123 & \( 2.466\)e\(+00\) & \( 3.748\)e\(-01\) & \( 1.776\)e\(+00\) \\
      cvxqp3\_s &   18 &     50 & \( 8.394\)e\(-02\) & \( 4.090\)e\(-02\) & \( 2.741\)e\(-02\) \\
      cvxqp3\_m &   19 &     81 & \( 4.505\)e\(-01\) & \( 1.737\)e\(-01\) & \( 2.369\)e\(-01\) \\
      cvxqp3\_l &   20 &    118 & \( 4.272\)e\(+00\) & \( 1.758\)e\(+00\) & \( 2.142\)e\(+00\) \\
       gouldqp3 &   10 &     12 & \( 5.600\)e\(-01\) & \( 4.870\)e\(-01\) & \( 5.158\)e\(-02\) \\
       gouldqp2 &    9 &      5 & \( 8.700\)e\(-02\) & \( 4.047\)e\(-02\) & \( 3.243\)e\(-02\) \\
       mosarqp1 &   17 &     39 & \( 6.036\)e\(-01\) & \( 2.265\)e\(-01\) & \( 3.166\)e\(-01\) \\
       mosarqp2 &   16 &     60 & \( 4.167\)e\(-01\) & \( 1.646\)e\(-01\) & \( 2.168\)e\(-01\) \\
         stcqp1 &   15 &    109 & \( 5.648\)e\(+00\) & \( 3.556\)e\(+00\) & \( 1.895\)e\(+00\) \\
         stcqp2 &   16 &    180 & \( 1.827\)e\(+00\) & \( 1.730\)e\(-01\) & \( 1.510\)e\(+00\) \\
\hline
\end{tabular}
\end{center}
}
\par \vspace{3pt}
\caption{Results for $K_{3.5}$ with solver CP-MINRES. Times are in seconds.\label{tab:res-k3.5-minres}}
\end{table}

We do not show the details for CP-DQGMRES(2) and CP-GMRES(2),
because they do not add much to the discussion. We summarize the results as follows: CP-DQGMRES(2) results in the same number of
PDCO and CP-Krylov iterations as CP-MINRES, as expected, and the corresponding times are comparable with those of MINRES.
CP-GMRES(2) results in an increase in the number of CP-Krylov iterations as compared with MINRES.
Whereas CP-DQGMRES with \(\ell > 2\) may be viewed as CP-MINRES with a form of partial reorthogonalization, setting \(\ell = 4\) did not yield any improvement on the symmetric formulations.

The results with CP-DQGMRES($\ell$) and CP-GMRES($\ell$) on $K_{3p}$ are not favorable.
In general, the unsymmetric CP-Krylov solvers on $K_{3p}$ are much less efficient than the symmetric ones on $K_2$ and $K_{3.5}$.
For example, with $\ell = 500$ the number of CPKrylov iterations is much larger than in the symmetric case and there are some problems
where CP-DQGMRES($\ell$) and CP-GMRES($\ell$) cannot always satisfy the stopping criterion. In these cases, they halt because a maximum
number of CP-Krylov iterations equal to $2n$ is achieved, thus preventing PDCO from computing the optimal solution by its maximum number
of iterations, which is set as $\min \{ \max \{ 30, n \}, \, 50 \}$. Among the possible reasons, we mention the non-normality of $K_{3p}$ and the
choice of the \((1,1)\) block $G$ of the preconditioner. A similar behavior has been observed by using the MATLAB
function \texttt{gmres} with the constraint preconditioner.
% \smarttodo{Actually, a better implementation (e.g., Householder orthogonalization)
% would reduce the number of inner iterations,
% as shown by the comparison with MATLAB GMRES, but the number of inner iterations is still too large.}
% I think an appropriate formulation will pay off much more than a more sophisticated implementation.
A more efficient choice of $G$ and a better formulation than \(K_{3p}\) are the subject of further investigation.

We ran all our tests a second time with iterative semi-refinement~\eqref{eq:proj-2-isr} activated, but did not observe any difference
in the number of inner or outer iterations.

\section{Discussion\label{sec:discussion}}

We extended the approach of~\cite{gould-orban-rees-2014} to saddle-point systems with regularization and provided principles from
which to derive constrained-preconditioned iterative methods.
The resulting methods are conceptually equivalent to standard iterative methods applied to a reduced system
in a way that preserves their properties, including quantities that increase or decrease monotonically at each iteration.
Specifically, we discussed constraint-preconditioned versions of the CG-Lanczos, MINRES, SYMMLQ, GMRES($\ell$)
and DQGMRES(\(\ell\)) methods, and showed that they preserve the properties of the corresponding standard methods in a suitable reduced Krylov space.
We illustrated our approach on methods based on the Lanczos and Arnoldi processes, but it applies equally to other processes,
including those of~\cite{golub-kahan-1965}, \cite{saunders-simon-yip-1988}, and the unsymmetric \cite{lanczos-1952} bi-orthogonalization process.
We implemented our constraint-preconditioned methods in a MATLAB library named \texttt{cpkrylov} that provides a basis for the development of more sophisticated numerical software.

An open question related to constraint preconditioners concerns the best way to reduce their computational cost. Inexact constraint preconditioners
have been developed and analyzed, based on approximations of the Schur complement of the leading block of the constraint preconditioner or on other
approximations~\citep{luksan-vlcek-1998, perugia-simoncini-2000, durazzi-ruggiero-2003,bergamaschi-gondzio-venturin-zilli-2007, sesana-simoncini-2013}.
Preconditioner updating techniques, producing inexact and exact constraint preconditioners, have been also proposed in order to reduce the cost of solving sequences of saddle-point systems \citep{bellavia-desimone-diserafino-morini-2015, bellavia-desimone-diserafino-morini-2016,
fisher-gratton-gurol-tremolet-vasseur-2016, bergamaschi-desimone-diserafino-martinez-2017}.
It must be noted, however, that the inexact constraint preconditioners considered so far generally do not produce preconditioned vectors lying in
the nullspace of the matrix \(N\) defined in~\eqref{eq:def-3x3}, which is a key issue to obtain CP-preconditioned methods
for~\eqref{eq:rsp} equivalent to suitably preconditioned Krylov methods for~\eqref{eq:reduced-sys}.
On the other hand, inexact preconditioners have proven effective
in reducing the computational time for the solution of large-scale saddle-point systems. A further possibility for lowering
the cost of constraint preconditioners is to apply them inexactly using an iterative method. Of course, preserving the property of obtaining
preconditioned vectors lying in the nullspace of \(N\) is a major issue. To the best of our knowledge, this approach has not been yet addressed
in the literature.

Finally, it is worth investigating the choice of the \((1,1)\) block of the constraint preconditioner when solving non-normal saddle-point systems.

\section*{Acknowledgments}
We thank two anonymous reviewers for constructive comments that helped us greatly improve the numerical experiments section.

\small
\bibliographystyle{abbrvnat}
\bibliography{cpkrylov}

\normalsize
\clearpage
\appendix

\section{Standard Lanczos and Arnoldi Processes\label{sec:processes}}

For reference we state the preconditioned Lanczos process and the full-space Lanczos process for~\eqref{eq:rsp} with preconditioner~\eqref{eq:cp}.
We also state the projected and full-space Arnoldi processes.

%\begin{algorithm}[ht!]
%    \caption{\label{alg:lanczos}%
%        Lanczos Process for $Ax = b$%
%    }
%    \begin{algorithmic}[1]
%%      \Require $A$, $b$, $x_0$
%      \State choose $x_0$
%      \State $v_0 = 0$, $k = 1$
%      \State $r_0 = b - A x_0$
%      \State $v_1 = r_0$
%      \State $\beta_1 = (v_1^T r_0)^{\tfrac{1}{2}}$
%      \If {$\beta_1 \ne 0$}
%          \State $v_1 = v_1 / \beta_1$ \Comment{$\|v_1\|_2 = 1$}
%      \EndIf
%      \State $k=1$
%      \While {$\beta_k \ne 0$}
%          \State $u_k = A v_k$
%          \State $\alpha_k = u_k^T v_k$
%          \State $v_{k+1} = u_k - \alpha_k v_k - \beta_k v_{k-1}$
%          \State $\beta_{k+1} = (v_{k+1}^T u_k)^{\tfrac{1}{2}}$
%          \If {$\beta_{k+1} \ne 0$}
%               \State $v_{k+1} = v_{k+1} / \beta_{k+1}$ \Comment{$\| v_{k+1} \|_2 = 1$}
%          \EndIf
%          \State $k = k+1$
%      \EndWhile
%    \end{algorithmic}
% \end{algorithm}
%%
%\vskip -10pt
\begin{algorithm}[hb!]
    \caption{\label{alg:prec-lanczos}%
        Lanczos Process for $Ax = b$ with Preconditioner $J = J^T \succ 0$%
    }
    \begin{algorithmic}[1]
%      \Require $A$, $P = P^T$, $b$, $x_0$ % \Comment{typically, $P \succ 0$}
      \State choose $x_0$
      \State $v_0 = 0$
      \State $r_0 = b - A x_0$
      \State solve $J v_1 = r_0$ % for $v_1$
      \State $\beta_1 = (v_1^T r_0)^{\tfrac{1}{2}}$
       \If {$\beta_1 \ne 0$}
          \State $v_1 = v_1 / \beta_1$ \Comment{$\|v_1\|_J = 1$}
      \EndIf
      \State $k=1$
      \While {$\beta_k \ne 0$}
          \State $u_k = A v_k$
          \State $\alpha_k = u_k^T v_k$
          \State solve $J v_{k+1} = u_k$ % for $v_{k+1}$
          \State $v_{k+1} = v_{k+1} - \alpha_k v_k - \beta_k v_{k-1}$
          \State $\beta_{k+1} = (v_{k+1}^T u_k)^{\tfrac{1}{2}}$
          \If {$\beta_{k+1} \ne 0$}
              \State $v_{k+1} = v_{k+1} / \beta_k$ \Comment{$\| v_{k+1} \|_J = 1$}
          \EndIf
          \State $k = k+1$
      \EndWhile
    \end{algorithmic}
\end{algorithm}

\begin{algorithm}[ht]
    \caption{\label{alg:cp-lanczos-full-space}%
        Full-Space Lanczos Process for~\eqref{eq:rsp} with Preconditioner~\eqref{eq:cp}%
    }
    \begin{algorithmic}[1]
%      \Require $A$, $B$, $C = C^T$, $G = G^T$, $b$, $[x_0\,;\, y_0]$ such that~\eqref{eq:2ndeq0_x0y0} holds
      \State choose $[x_0\,;\, y_0]$ such that~$B x_0 - C y_0 = 0$
      \State initialize
        \[
          \begin{bmatrix}
            v_{0,x} \\ v_{0,y}
          \end{bmatrix}
          =
          \begin{bmatrix}
            0 \\ 0
          \end{bmatrix}
        \]
      \State set \Comment{$B x_0 - C y_0 = 0 \; \Rightarrow \; r_{0,y} = 0$}
        \[
          \begin{bmatrix}
            r_{0,x} \\ r_{0,y}
          \end{bmatrix}
          =
          \begin{bmatrix}
            b \\ 0
          \end{bmatrix}
          -
          \begin{bmatrix}
            A & \phantom{-}B^T \\
            B & -C\phantom{^T}
          \end{bmatrix}
          \begin{bmatrix}
            x_0 \\ y_0
          \end{bmatrix}
          =
          \begin{bmatrix}
            b - A x_0 - B^T y_0 \\ 0
          \end{bmatrix}
        \]
      \State \label{line:v1-full-space}%
       obtain $v_1$ as the solution of
        \[
          \begin{bmatrix}
            G & \phantom{-}B^T \\
            B & -C\phantom{^T}
          \end{bmatrix}
          \begin{bmatrix}
            v_{1,x} \\ v_{1,y}
          \end{bmatrix}
          =
          \begin{bmatrix}
            r_{0,x} \\ r_{0,y}
          \end{bmatrix}
        \]
      \State \label{alg:cplfs:beta1}%
        $\beta_1 = (v_1^T r_0)^{\tfrac{1}{2}} = (v_{1,x}^T r_{0,x})^{\tfrac{1}{2}}$
%      \IfOneLine {$\beta_1 \ne 0$} \ $v_1 = v_1 / \beta_1$
      \If {$\beta_1 \ne 0$}
           \State $v_1 = v_1 / \beta_1$
      \EndIf
      \State $k = 1$
      \While {$\beta_k \ne 0$}
          \State compute
            \[
              \begin{bmatrix}
                u_{k,x} \\ u_{k,y}
              \end{bmatrix}
              =
              \begin{bmatrix}
                A & \phantom{-}B^T \\
                B & -C\phantom{^T}
              \end{bmatrix}
              \begin{bmatrix}
                v_{k,x} \\ v_{k,y}
              \end{bmatrix}
            \]
          \State \label{alg:cplfs:alphak}%
            $\alpha_k = u_k^T v_k = v_{k,x}^T A v_{k,x} + 2 v_{k,x}^T B^T v_{k,y} - v_{k,y}^T C v_{k,y}$
          \State \label{alg:cplfs:vk+1}%
            obtain \(v_{k+1}\) as the solution of
            \[
              \begin{bmatrix}
                G & \phantom{-}B^T \\
                B & -C\phantom{^T}
              \end{bmatrix}
              \begin{bmatrix}
                v_{k+1,x} \\ v_{k+1,y}
              \end{bmatrix}
              =
              \begin{bmatrix}
                u_{k,x} \\ u_{k,y}
              \end{bmatrix}
            \]
          \State \label{alg:cplfs:update-vk+1}%
            $v_{k+1} = v_{k+1} - \alpha_k v_k - \beta_k v_{k-1}$
          \State $\beta_{k+1} = (v_{k+1}^T u_k)^{\tfrac{1}{2}} = (v_{k+1,x}^T u_{k,x} + v_{k+1,y}^T u_{k,y})^{\tfrac{1}{2}}$
%          \IfOneLine {$\beta_{k+1} \ne 0$} \ $v_{k+1} = v_{k+1} / \beta_{k+1}$
          \If {$\beta_{k+1} \ne 0$}
               \State $v_{k+1} = v_{k+1} / \beta_{k+1}$
          \EndIf
          \State $k = k+1$
      \EndWhile
    \end{algorithmic}
\end{algorithm}

\begin{algorithm}[ht]
   \caption{Projected Arnoldi Process \label{alg:proj-arnoldi}} % for~\eqref{eq:rsp-2}
   \begin{algorithmic}[1]
     \State choose $[x_0\,;\, w_0]$ such that $B x_0 + E w_0 = 0$
     \State $v_{0,x} = 0$, \  $v_{0,w} = -w_0$
     \State $u_{0,x} = b - A x_0$, $u_{0,w} = -F^{-1} w_0$
     \State $[\bar{u}_{1,x} \,; \,\bar{u}_{1,w} \,;\, \bar{z}_1] \leftarrow$ solution of~\eqref{eq:proj} with right-hand side $[u_{0,x}\,;\, u_{0,w}\,;\, 0]$
               \label{alg:a-line-proj-v1}
     \State $v_{1,x} = \bar{u}_{1,x}$, \ $v_{1,w} = \bar{u}_{1,w}$
               \Comment{$v_1 = P_G \, u_0$}
     \State $h_{1,0} = (v_{1,x}^T u_{0,x} + v_{1,w}^T u_{0,w})^{\tfrac{1}{2}}$
     \If {$h_{1,0} \ne 0$}
         \State $v_{1,x} = v_{1,x} / h_{1,0}$,  \ $v_{1,w} = v_{1,w} / h_{1,0}$  \label{alg:a-line-scale-v1}
     \EndIf
     \State $k = 1$
     \While {$h_{k,k-1}  \ne 0$}
         \State $u_{k,x} = A v_{k,x}$, \ $u_{k,w} = F^{-1} v_{k,w}$
         \State $[\bar{u}_{k+1,x}\,;\, \bar{u}_{k+1,w}\,;\, \bar{z}_{k+1}] \leftarrow$ solution of~\eqref{eq:proj} with right-hand side $[u_{k,x}\,;\, u_{k,w}\,;\, 0]$
                    \label{alg:a-line-proj-vk+1}
         \State $v_{k+1,x} = \bar{u}_{k+1,x}$, \ $v_{k+1,w} = \bar{u}_{k+1,w}$
                    \Comment{$v_{k+1} = P_G \, u_k$}
         \For {$i = 1,\ldots,k$}
             \State $h_{i,k} = v_{i,x}^T u_{k,x} + v_{i,w}^T u_{k,w}$
             \State $v_{k+1,x} = v_{k+1,x} - h_{i,k} v_{i,x}$,  \ $v_{k+1,w} = v_{k+1,w} - h_{i,k} v_{i,w}$
         \EndFor
         \State $h_{k+1,k} = (v_{k+1,x}^T u_{k,x} + v_{k+1,w}^T u_{k,w})^{\tfrac{1}{2}}$
         \If {$h_{k+1,k} \ne 0$}
             \State $v_{k+1,x} = v_{k+1,x} / h_{k+1,k}$, \ $v_{k+1,w} = v_{k+1,w} / h_{k+1,k}$
                       \label{alg:a-line-scale-vk+1}
         \EndIf
         \State $k = k+1$
     \EndWhile
   \end{algorithmic}
\end{algorithm}

\begin{algorithm}[ht]
   \caption{\label{alg:cp-arnoldi-full-space}%
       Full-Space Arnoldi Process for~\eqref{eq:rsp} with Preconditioner~\eqref{eq:cp}%
   }
   \begin{algorithmic}[1]
%     \Require $A$, $B$, $C = C^T$, $G = G^T$, $b$, $[x_0\,;\, y_0]$ such that~\eqref{eq:2ndeq0_x0y0} holds
     \State choose $[x_0\,;\, y_0]$ such that $B x_0 - C y_0 = 0$
     \State initialize
       \[
         \begin{bmatrix}
           v_{0,x} \\ v_{0,y}
         \end{bmatrix}
         =
         \begin{bmatrix}
           0 \\ 0
         \end{bmatrix}
       \]
     \State set \Comment{$B x_0 - C y_0 = 0 \; \Rightarrow \; r_{0,y} = 0$}
       \[
         \begin{bmatrix}
           r_{0,x} \\ r_{0,y}
         \end{bmatrix}
         =
         \begin{bmatrix}
           b \\ 0
         \end{bmatrix}
         -
         \begin{bmatrix}
           A & \phantom{-}B^T \\
           B & -C\phantom{^T}
         \end{bmatrix}
         \begin{bmatrix}
           x_0 \\ y_0
         \end{bmatrix}
         =
         \begin{bmatrix}
           b - A x_0 - B^T y_0 \\ 0
         \end{bmatrix}
       \]
     \State obtain $v_1$ as the solution of
       \[
         \begin{bmatrix}
           G & \phantom{-}B^T \\
           B & -C\phantom{^T}
         \end{bmatrix}
         \begin{bmatrix}
           v_{1,x} \\ v_{1,y}
         \end{bmatrix}
         =
         \begin{bmatrix}
           r_{0,x} \\ r_{0,y}
         \end{bmatrix}
       \]
     \State $h_{1,0} = (v_1^T r_0)^{\tfrac{1}{2}} = (v_{1,x}^T r_{0,x})^{\tfrac{1}{2}}$
%     \IfOneLine {$h_{1,0} \ne 0$} \ $v_1 = v_1 / \beta_1$
     \If {$h_{1,0} \ne 0$}
         \State $v_1 = v_1 / \beta_1$
     \EndIf
     \State $k = 1$
     \While {$h_{k,k-1} \ne 0$}
         \State compute
           \[
             \begin{bmatrix}
               u_{k,x} \\ u_{k,y}
             \end{bmatrix}
             =
             \begin{bmatrix}
               A & \phantom{-}B^T \\
               B & -C\phantom{^T}
             \end{bmatrix}
             \begin{bmatrix}
               v_{k,x} \\ v_{k,y}
             \end{bmatrix}
           \]
         \State obtain \(v_{k+1}\) as the solution of
           \[
             \begin{bmatrix}
               G & \phantom{-}B^T \\
               B & -C\phantom{^T}
             \end{bmatrix}
             \begin{bmatrix}
               v_{k+1,x} \\ v_{k+1,y}
             \end{bmatrix}
             =
             \begin{bmatrix}
               u_{k,x} \\ u_{k,y}
             \end{bmatrix}
           \]
         \For {\(i = 1, \dots, k\)}
           \State $h_{i,k} = v_i^T u_k = v_{i,x}^T A v_{k,x} + 2 v_{i,x}^T B^T v_{k,y} - v_{i,y}^T C v_{k,y}$
           \State $v_{k+1} = v_{k+1} - h_{i,k} v_i$
         \EndFor
         \State $h_{k+1,k} = (v_{k+1}^T u_k)^{\tfrac{1}{2}} = (v_{k+1,x}^T u_{k,x} + v_{k+1,y}^T u_{k,y})^{\tfrac{1}{2}}$
%         \IfOneLine {$h_{k+1,k} \neq 0$} \ $v_{k+1} = v_{k+1} / h_{k+1,k}$
         \If {$h_{k+1,k} \neq 0$}
             \State $v_{k+1} = v_{k+1} / h_{k+1,k}$
         \EndIf
         \State $k = k+1$
     \EndWhile
   \end{algorithmic}
\end{algorithm}

\clearpage

%\tableofcontents
% \listoftodos\relax
% \listofchanges

%\smarttodo[color=white,linecolor=gray]{...}

\end{document}